\documentclass[12pt]{amsart}
\usepackage{amsmath}
\usepackage{amsfonts}
\usepackage{amssymb}
\usepackage{mathtools}
\usepackage{bbding}
\usepackage{stmaryrd}
\usepackage{txfonts}
\usepackage{graphicx}
\usepackage{epsfig}
\usepackage{xypic}
\usepackage{tikz}
\usepackage{longtable}
\usepackage[all]{xy}
\usepackage{pgflibraryarrows}
\usepackage{pgflibrarysnakes}
\usepackage[shortlabels]{enumitem}
\usepackage{ifpdf}
\ifpdf
\usepackage[colorlinks,final,backref=page,hyperindex]{hyperref}
\else
\usepackage[colorlinks,final,backref=page,hyperindex,hypertex]{hyperref}
\fi
\usepackage{graphicx}
\usepackage{epstopdf}
\usepackage{epsfig}
\usepackage{pdfsync}    

\usepackage{xspace}

\newtheorem{theorem}{Theorem}[section]
\newtheorem{lemma}[theorem]{Lemma}
\newtheorem{corollary}[theorem]{Corollary}

\newtheorem{proposition}[theorem]{Proposition}

\theoremstyle{definition}
\newtheorem{definition}[theorem]{Definition}
\newtheorem{remark}[theorem]{Remark}
\newtheorem{example}[theorem]{Example}

\newcommand{\nc}{\newcommand}
\newcommand{\delete}[1]{}

\delete{

}

\def\bc{\begin{center}}
	\def\ec{\end{center}}

\nc{\tred}[1]{\textcolor{red}{#1}}
\nc{\tblue}[1]{\textcolor{blue}{#1}} \nc{\tgreen}[1]{\textcolor{green}{#1}} \nc{\tpurple}[1]{\textcolor{purple}{#1}} \nc{\btred}[1]{\textcolor{red}{\bf #1}} \nc{\btblue}[1]{\textcolor{blue}{\bf #1}} \nc{\btgreen}[1]{\textcolor{green}{\bf #1}} \nc{\btpurple}[1]{\textcolor{purple}{\bf #1}}

\newcommand{\efootnote}[1]{}

\nc{\mlabel}[1]{\label{#1}}  
\nc{\mcite}[1]{\cite{#1}}  
\nc{\mref}[1]{\ref{#1}}  
\nc{\meqref}[1]{\eqref{#1}}  
\nc{\mbibitem}[1]{\bibitem{#1}} 

\delete{
	\nc{\mlabel}[1]{\label{#1}  
		{\hfill \hspace{1cm}{\bf{{\ }\hfill(#1)}}}}
	\nc{\mcite}[1]{\cite{#1}{{\bf{{\ }(#1)}}}}  
	\nc{\mref}[1]{\ref{#1}{{\bf{{\ }(#1)}}}}  
	\nc{\meqref}[1]{\eqref{#1}{{\bf{{\ }(#1)}}}}  
	\nc{\mbibitem}[1]{\bibitem[\bf #1]{#1}} 
}

\renewcommand\geq{\geqslant}
\renewcommand\leq{\leqslant}

\renewcommand\bar[1]{\overline{#1}}

\nc{\name}[1]{{\bf #1}}

\nc{\tforall}{\quad \text{ for all }}

\nc{\mre}{\text{Re}\,}

\nc{\mim}{\text{im}\,}
\nc{\nz}{\varepsilon}
\nc{\Id}{\mathrm{Id}}

\nc{\DO}{\text{DO}}
\nc{\IDO}{\text{IDO}}

\nc{\IEO}{\mathrm{IEO}}

\nc{\mnoindent}{\smallskip\noindent}

\nc{\lnvkv}{\triangleleft}

\nc{\bin}[2]{ (_{\stackrel{\scs{#1}}{\scs{#2}}})}  
\nc{\binc}[2]{ \left (\!\! \begin{array}{c} \scs{#1}\\
		\scs{#2} \end{array}\!\! \right )}  
\nc{\bincc}[2]{  \left ( {\scs{#1} \atop
		\vspace{-1cm}\scs{#2}} \right )}  
\nc{\bs}{\bar{S}} \nc{\cosum}{\sqsubset} \nc{\la}{\longrightarrow} \nc{\rar}{\rightarrow} \nc{\dar}{\downarrow} \nc{\dprod}{**} \nc{\dap}[1]{\downarrow \rlap{$\scriptstyle{#1}$}} \nc{\md}[1]{\bar{#1}} \nc{\uap}[1]{\uparrow \rlap{$\scriptstyle{#1}$}} \nc{\defeq}{\stackrel{\rm def}{=}} \nc{\disp}[1]{\displaystyle{#1}} \nc{\dotcup}{\ \displaystyle{\bigcup^\bullet}\ } \nc{\gzeta}{\bar{\zeta}} \nc{\hcm}{\ \hat{,}\ } \nc{\hts}{\hat{\otimes}} \nc{\barot}{{\otimes}} \nc{\free}[1]{\bar{#1}} \nc{\uni}[1]{\tilde{#1}} \nc{\hcirc}{\hat{\circ}} \nc{\leng}{\ell} \nc{\lleft}{[} \nc{\lright}{]} \nc{\lc}{\lfloor} \nc{\rc}{\rfloor}
\nc{\curlyl}{\left \{ \begin{array}{c} {} \\ {} \end{array}
	\right.  \!\!\!\!\!\!\!}
\nc{\curlyr}{ \!\!\!\!\!\!\!
	\left. \begin{array}{c} {} \\ {} \end{array}
	\right \} }
\nc{\longmid}{\left | \begin{array}{c} {} \\ {} \end{array}
	\right. \!\!\!\!\!\!\!}
\nc{\onetree}{\bullet} \nc{\ora}[1]{\stackrel{#1}{\rar}}
\nc{\ola}[1]{\stackrel{#1}{\la}}
\nc{\ot}{\otimes} \nc{\mot}{{{\boxtimes\,}}} \nc{\otm}{\overline{\boxtimes}} \nc{\sprod}{\bullet} \nc{\scs}[1]{\scriptstyle{#1}} \nc{\mrm}[1]{{\rm #1}} \nc{\msum}{\sum\limits}
\nc{\margin}[1]{\marginpar{\rm #1}}   
\nc{\dirlim}{\displaystyle{\lim_{\longrightarrow}}\,} \nc{\invlim}{\displaystyle{\lim_{\longleftarrow}}\,} \nc{\mvp}{\vspace{0.3cm}} \nc{\tk}{^{(k)}} \nc{\tp}{^\prime} \nc{\ttp}{^{\prime\prime}} \nc{\svp}{\vspace{2cm}} \nc{\vp}{\vspace{8cm}} \nc{\proofbegin}{\noindent{\bf Proof: }}
\nc{\proofend}{$\blacksquare$ \vspace{0.3cm}}
\nc{\modg}[1]{\!<\!\!{#1}\!\!>}
\nc{\intg}[1]{F_C(#1)} \nc{\lmodg}{\!<\!\!} \nc{\rmodg}{\!\!>\!} \nc{\cpi}{\widehat{\Pi}}
\nc{\sha}{{\mbox{\cyr X}}}  
\nc{\shap}{{\mbox{\cyrs X}}} 
\nc{\shpr}{\diamond}    
\nc{\shp}{\ast} \nc{\shplus}{\shpr^+}
\nc{\shprc}{\shpr_c}    
\nc{\msh}{\ast} \nc{\zprod}{m_0} \nc{\oprod}{m_1} \nc{\vep}{\varepsilon} \nc{\labs}{\mid\!} \nc{\rabs}{\!\mid}
\nc{\astarrow}{\overset{\raisebox{-3pt}{$\ast$}}{\rightarrow}}

\nc{\sqsym}{Stirling quasisymmetric function\xspace}
\nc{\sqsyms}{Stirling quasisymmetric functions\xspace}

\nc{\EEsym}{\mathbb{E}sym}
\nc{\Sym}{\mrm{Sym}}
\nc{\NSym}{\mrm{NSym}}
\nc{\QSym}{\mrm{QSym}}
\nc{\RQSym}{\mrm{RQSym}}
\nc{\RenQSym}{\mrm{WCQSym}}	
\nc{\DQSym}{\mrm{DQSym}}
\nc{\WDQSym}{\mrm{WDQSym}}
\nc{\DLQSym}{\mrm{DLQSym}}
\nc{\ZQSym}{\mrm{ZQSym}}
\nc{\Ensym}{\mrm{ENSym}}
\nc{\Wcsym}{\mrm{WCSym}}
\nc{\LWQSym}{\mrm{LWQSym}}
\nc{\LWCQSym}{\mrm{\mathrm{LWQSym}}}
\nc{\Wcqsym}{\mrm{QSym}_{\widetilde{\mathbb{N}}}}
\nc{\Syms}{symmetric functions\xspace}
\nc{\eqsym}{extended quasisymmetric function\xspace}
\nc{\eqsyms}{extended quasisymmetric functions\xspace}
\nc{\Eqsyms}{Extended Quasisymmetric functions\xspace}
\nc{\Esyms}{Extended symmetric functions\xspace}
\nc{\sgqsym}{quasisymmetric function with semigroup exponents\xspace}
\nc{\sgqsyms}{quasisymmetric functions with semigroup exponents\xspace}
\nc{\Sgqsyms}{Quasisymmetric functions with semigroup exponents\xspace}
\nc{\SGQSYM}{\mrm{SGQSYM}}
\nc{\emzv}{extended multiple zeta value}
\nc{\emzvs}{extended multiple zeta values}
\nc{\sgfps}{formal power series with semigroup exponent\xspace}
\nc{\NSymg}{\mathrm{NSym}_\gp}
\nc{\zqsym}{zeta-quasisymmetric }
\nc{\gslwqsym}{Stirling left weak quasisymmetric function\xspace}
\nc{\gslwqsyms}{Stirling left weak quasisymmetric functions\xspace}
\nc{\ulwb}{upper-left weak bicomposition\xspace}
\nc{\ulwbs}{upper-left weak bicompositions\xspace}

\nc{\parr}{\rm Par}
\nc{\wpar}{\rm WPar}
\nc{\wcomp}{\large{\VDash}}
\nc{\Ker}{\ker}

\nc{\dth}{d} \nc{\mmbox}[1]{\mbox{\ #1\ }} \nc{\fp}{\mrm{FP}} \nc{\rchar}{\mrm{char}} \nc{\Fil}{\mrm{Fil}} \nc{\Mor}{Mor\xspace} \nc{\gmzvs}{gMZV\xspace} \nc{\gmzv}{gMZV\xspace} \nc{\mzv}{MZV\xspace} \nc{\mzvs}{MZVs\xspace}
\nc{\MZV}{\mathrm{MZV}}
\nc{\Hom}{\mrm{Hom}} \nc{\id}{\mrm{id}} \nc{\im}{\mrm{im}} \nc{\incl}{\mrm{incl}}  \nc{\mchar}{\rm char}

\nc{\Alg}{\mathbf{Alg}} \nc{\Bax}{\mathbf{Bax}} \nc{\bff}{\mathbf f} \nc{\bfk}{{\bf k}} \nc{\bfone}{{\bf 1}} \nc{\bfx}{\mathbf x} \nc{\bfy}{\mathbf y}
\nc{\base}[1]{\bfone^{\otimes ({#1}+1)}} 
\nc{\Cat}{\mathbf{Cat}} \delete{}
\nc{\detail}{\marginpar{\bf More detail}
	\noindent{\bf Need more detail!}
	\svp}
\nc{\Int}{\mathbf{Int}} \nc{\Mon}{\mathbf{Mon}}
\nc{\rbtm}{{shuffle }} \nc{\rbto}{{Rota-Baxter }} \nc{\remarks}{\noindent{\bf Remarks: }} \nc{\Rings}{\mathbf{Rings}} \nc{\Sets}{\mathbf{Sets}}
\nc{\balpha}{\mathbf{\alpha}}

\nc{\BA}{{\mathbb A}} \nc{\CC}{{\mathbb C}} \nc{\DD}{{\mathbb D}} \nc{\EE}{{\mathbb E}} \nc{\FF}{{\mathbb F}} \nc{\GG}{{\mathbb G}} \nc{\HH}{{\mathbb H}} \nc{\LL}{{\mathbb L}} \nc{\NN}{{\mathbb N}} \nc{\KK}{{\mathbb K}} \nc{\PP}{{\mathbb P}} \nc{\QQ}{{\mathbb Q}} \nc{\RR}{{\mathbb R}} \nc{\TT}{{\mathbb T}} \nc{\VV}{{\mathbb V}} \nc{\ZZ}{{\mathbb Z}}

\nc{\cala}{{\mathcal A}} \nc{\calc}{{\mathcal C}} \nc{\cald}{{\mathcal D}} \nc{\cale}{{\mathcal E}} \nc{\calf}{{\mathcal F}} \nc{\calg}{{\mathcal G}} \nc{\calh}{{\mathcal H}} \nc{\cali}{{\mathcal I}} \nc{\call}{{\mathcal L}} \nc{\calm}{{\mathcal M}} \nc{\caln}{{\mathcal N}} \nc{\calo}{{\mathcal O}} \nc{\calp}{{\mathcal P}} \nc{\calr}{{\mathcal R}} \nc{\cals}{{\mathcal S}} \nc{\calt}{{\mathcal T}} \nc{\calw}{{\mathcal W}} \nc{\calk}{{\mathcal K}} \nc{\calx}{{\mathcal X}}
\nc{\calz}{{\mathcal Z}}

\nc{\fraka}{{\mathfrak a}} \nc{\frakA}{{\mathfrak A}} \nc{\frakb}{{\mathfrak b}} \nc{\frakB}{{\mathfrak B}}
\nc{\frakc}{{\mathfrak c}}  \nc{\frakD}{{\mathfrak D}}
\nc{\frakH}{{\mathfrak H}}
\nc{\frakh}{{\mathfrak h}} \nc{\frakM}{{\mathfrak M}}
\nc{\frakO}{{\mathfrak O}}
\nc{\frakE}{{\mathfrak E}}
\nc{\bfrakM}{\overline{\frakM}} \nc{\frakm}{{\mathfrak m}} \nc{\frakP}{{\mathfrak P}} \nc{\frakN}{{\mathfrak N}} \nc{\frakp}{{\mathfrak p}} \nc{\frakS}{{\mathfrak S}}
\nc{\frakk}{{\mathfrak k}}
\nc{\frakx}{{\mathfrak x}}
\nc{\frakl}{{\mathfrak l}} \nc{\ox}{\bar{\frakx}} \nc{\frakX}{{\mathfrak X}} \nc{\fraky}{{\mathfrak y}} \nc\dop{\delta}
\nc{\Reduce}{{\rm Red}}

\font\cyr=wncyr10 \font\cyrs=wncyr7
\nc{\redt}[1]{\textcolor{red}{#1}}
\nc{\li}[1]{\textcolor{red}{#1}}
\nc{\lir}[1]{\textcolor{red}{Li:#1}}
\nc{\ap}[1]{\textcolor{blue}{#1}}
\nc{\apr}[1]{\textcolor{blue}{AP:#1}}

\nc{\Supp}{\mathrm{Supp}}
\nc{\Fix}{{\mathrm{Fix}}}

\nc{\SD}{{\text{SD}}}

\nc{\Ideal}{{\mathrm{Ideal}}}

\nc{\diam}{{\mathrm{diam}}}
\nc{\gr}{{\mathrm{gr}}}
\nc{\D}{{\mathrm{D}}}

\nc{\rb}{integral\xspace}
\nc{\Rb}{Integral\xspace}
\nc{\mrb}{Rota-Baxter\xspace}

\nc{\wvec}[2]{{\scriptsize{\Big [ \!\!\begin{array}{c} #1 \\ #2 \end{array} \!\! \Big ]}}}

\nc{\bwvec}[2]{\Big(\wvec{#1}{#2}\Big)}

\nc{\jwvec}[2]{{\scriptsize{\Big [ \!\!\begin{array}{cccccccccccccc} #1 \\ #2 \end{array} \!\! \Big ]}}}
\nc{\bjwvec}[2]{\Big(\jwvec{#1}{#2}\Big)}

\begin{document}
	
\title[Ideal-based Zero-divisor Graph of MV-algebras]{Ideal-based Zero-divisor Graph of MV-algebras}

\author{Aiping Gan}
\address{School of Mathematics and Statistics,
Jiangxi Normal University, Nanchang, Jiangxi 330022, P.R. China}
\email{ganaiping78@163.com}

\author{Huadong Su}
\address{School of Science,
	Beibu Gulf University, Qinzhou, Guangxi 535011, P.R. China}
\email{huadongsu@sohu.com}

\author{Yichuan Yang}
\address{Department of Mathematics,
	Beihang University, Beijing, 100191, China}
\email{ycyang@buaa.edu.cn(corresponding author)}

\hyphenpenalty=8000
	
\date{\today}
	
\begin{abstract} 
Let $(A, \oplus, *, 0)$ be an MV-algebra,  $(A, \odot, 0)$ be the associated commutative semigroup,
and $I$ be an ideal of $A$. Define the ideal-based zero-divisor graph
$\Gamma_{I}(A)$ of $A$ with respect to $I$ to be a simple graph with
the set of vertices
$V(\Gamma_{I}(A))=\{x\in A\backslash I ~|~ (\exists~ y\in A\backslash I) ~x\odot y\in I\},$
and  two distinct vertices $x$ and $y$ are joined by an edge if and only
if $x\odot y\in I$.
 We prove that
$\Gamma_{I}(A)$ is connected and its diameter is less than or equal to
$3$.  Also, some
relationship between the diameter (the girth) of $\Gamma_{I}(A)$ and the diameter (the girth) of the zero-divisor graph of $A/I$ are investigated.
And using the girth of  zero-divisor graphs (resp. ideal-based zero-divisor graphs) of MV-algebras, we classify all MV-algebras into $2~($resp. $3)$ types.
\end{abstract}
	
\subjclass[2010]{
	06D35, 
	 05C25, 
05C38  
05C40 
}
	
\keywords{MV-algebra; zero-divisor graph, ideal-based zero-divisor graph, girth, diameter}
	
\maketitle
	
\tableofcontents
	
\hyphenpenalty=8000 \setcounter{section}{0}
	
	
\allowdisplaybreaks

\section{Introduction and preliminary}
Let $R$ be a commutative ring with nonzero identity. The zero-divisor graphs of  $R$ was first introduced by Beck \cite{beck}, where he was mostly concerned in coloring. In his work all elements of the ring were vertices of the graph and two distinct elements $x$ and $y$ are adjacent if and only if $xy=0$.  
Anderson and Livingston
\cite{anderson2} 
have defined a  graph $\Gamma(R)$ of $R$ to be a simple graph whose vertices are the nonzero
zero-divisors of $R$, and two distinct vertices $x$ and $y$ are joined by an edge
if and only if $xy=0$. 
They completely characterized the star graphs
of the form $\Gamma(R)$ where $R$ is finite, proving that the star graphs $G$ that occurs as $\Gamma(R)$
are precisely those such that the cardinality $|G|$ of $G$ is a prime power.
 Since then
the interplay between algebraic properties of a ring $R$ and the graph theoretic properties of
the zero-divisor graph $\Gamma(R)$ has been studied extensively by many authors, see e.g.
 \cite{AMY,LW}.

 Redmond \cite{red} extended the zero-divisor graph of a commutative ring to
an ideal-based zero-divisor graph of a commutative ring.
Let $I$ be an ideal of a commutative ring $R$. The ideal-based zero-divisor graph (with respect to
$I$) of $R$, denoted by $\Gamma_{I}(R)$,
is a simple graph  with the set of vertices $$V(\Gamma_{I}(R))=\{x\in R\backslash I ~|~ (\exists ~y\in R\backslash I) ~xy\in I\}$$
where $R\backslash I=\{x\in R ~|~x\not\in I\}$, and
two distinct vertices $x$ and $y$ are joined by an edge if and only
if $xy\in I$. For more details on ideal-based zero-divisor graphs of  commutative rings, please see
\cite{jgs,maim,nik,red,smith}.
The (ideal-based) zero-divisor graph has also been extended to other algebraic structure,
such as zero-divisor graphs of semigroups  \cite{de2,de1}, ideal-based zero-divisor graphs of commutative semirings
 \cite{at1,at2}.

 Gan and Yang \cite{gan} introduced  the zero-divisor graphs of MV-algebras. 
Let $(A, \oplus, *, 0)$ be an MV-algebra, 
and $(A, \odot, 0)$ be the associated commutative semigroup. 
The zero-divisor graph $\Gamma(A)$ of  $A$ is defined to be a simple graph
whose vertices are the nonzero
zero-divisors of $A$, and two distinct vertices $x$ and $y$ are joined by an edge
if and only if $x\odot y=0$. 
They prove that the set of vertices of $\Gamma(A)$ is just $A\backslash\{0, 1\}$\cite[Proposition 2]{gan}. So $\Gamma(A)$
would be very large when the the cardinality of $A$ is large enough. To compress $\Gamma(A)$, 
we will focus on the ideal-based zero-divisor graphs of  MV-algebras in this paper.

The paper is organized as follows. In Subsections \ref{subsec:1.1} and \ref{subsec:2.2},
we review some necessary notions about MV-algebras and graphs.
In Section \ref{sec:3}, some properties about ideals of MV-algebras  are obtained. 
We prove that if an MV-algebra $A$  
is the direct product  of a family $\{A_{i}\}_{i\in \Omega}$ of MV-algebras, and
 $I_{i}$ be an ideal of $A_{i}$ for each $i\in \Omega$, then  the cartesian  product $\prod_{i\in \Omega}I_{i}$ of sets $I_{i}$  is an ideal of $A$, and $A/(\prod_{i\in \Omega}I_{i})\cong \prod_{i\in \Omega}(A_{i}/I_{i})$ (Theorem \ref{the:1000}). Moreover, when $\Omega$ is  finite, every ideal of $A$ can be expressed as the form: the cartesian  product $\prod_{i\in \Omega}I_{i}$, where $I_{i}$ is an ideal of $A_{i}$ for each $i\in \Omega$ (Proposition \ref{cor:10000}).
 In Section \ref{sec:4}, we investigate the girth of the zero-divisor graph of an MV-algebra, and obtain that the girth of $\Gamma(A)$ is equal to $3$ or $\infty$ (Theorem \ref{the:4.30}).
In Section \ref{sec:5}, we introduce and study the ideal-based zero-divisor graph of an MV-algebra.
Let  $I$ be an ideal of an MV-algebra $A$.
The ideal-based zero-divisor graph
$\Gamma_{I}(A)$ of $A$ with respect to $I$ is defined to be a simple graph with
the set of vertices
$$V(\Gamma_{I}(A))=\{x\in A\backslash I ~|~ (\exists~ y\in A\backslash I) ~x\odot y\in I\},$$
and  two distinct vertices $x$ and $y$ are joined by an edge
if and only if $x\odot y\in I$. It is clear that $\Gamma_{I}(A)=\Gamma(A)$ when $I=\{0\}$.
We prove that
$\Gamma_{I}(A)$ is connected and its diameter is less then or equal to
$3$ (Theorem \ref{the:1}). Also, some
relationship between the diameter of $\Gamma_{I}(A)$ and the diameter of the zero-divisor graph $\Gamma(A/I)$ of $A/I$ are investigated: when $I\neq \{0\}$, the diameter of $\Gamma_{I}(A)$  equals to $2$ if and only if the diameter of $\Gamma (A/I)$  equals to $1$ or $2$ (Theorem \ref{the:5}); and the diameter of $\Gamma_{I}(A)$  equals to $3$ if and only if the diameter of $\Gamma (A/I)$  is $3$ (Theorem \ref{the:6}). Finally, the girth of $\Gamma_{I}(A)$ is discussed, and using the girth of $\Gamma_{I}(A)$, we classify all MV-algebras into 3 types (Theorem \ref{the:5.20}).


\subsection{Some notions of MV-algebras}
\label{subsec:1.1}

\begin{definition} \label{def:1}
	\cite{cc,mu}~An \name{MV-algebra} is an algebra $(A, \oplus, *, 0)$  of type $(2, 1, 0)$ satisfying the following axioms: for all $x, y\in A$,
	\begin{enumerate}
		\item [$(M1)$]   $(A, \oplus, 0)$ is a commutative monoid;
		\item [$(M2)$]  $x^{**}=x$;
		\item [$(M3)$] $x\oplus 0^{*}=0^{*}$;
		\item [$(M4)$] $(x^{*}\oplus y)^{*}\oplus y=(y^{*}\oplus x)^{*}\oplus x$.
	\end{enumerate}
\end{definition}

As usual we shall denote an MV-algebra by its underlying carrier set.
Note that all axioms of MV-algebras are equations between terms,
MV-algebras form a variety.
So the notions of isomorphism, subalgebra, congruence, direct product
are just the particular cases of  the corresponding universal algebraic notions \cite{bu}.
In the rest of the paper, unless otherwise specified, we always assume that
$A=(A, \oplus, *, 0)$ is a \name{nontrivial  MV-algebra}, i.e, $A$ is an MV-algebra containing at least two elements.

Define the constant $1$ and the operation
$\odot$  on $A$ as:
	$1 =0^{*}$ and $ x\odot y=(x^{*}\oplus y^{*})^{*}$.
Then for all $x, y\in A$, the following 
properties hold \cite{mu}:
\begin{enumerate}
	\item [$\bullet$]    $(A, \odot, *, 1)$ is  an MV-algebra;
	\item [$\bullet$]  $\ast$ is an
	isomorphism between $(A, \oplus, *, 0)$ and  $(A, \odot, *, 1)$,
	\item [$\bullet$]  $1^{*}=0$;
	\item [$\bullet$] $ x\oplus y=(x^{*}\odot y^{*})^{*}$;
	\item [$\bullet$] $ x\oplus 1=1$;
	\item [$\bullet$] $x\oplus x^{*}=1$;
	\item [$\bullet$] $x\odot x^{*}=0$.
\end{enumerate}

\begin{example}\label{exa:1}
	\emph{\cite{cc,mu}}~
	Equip the real unit interval $[0, 1]$  with the
	operations $$x\oplus y = \min\{1, x+ y\} ~~\textrm{and}~~ x^{*}=1-x.$$
	Then $[0, 1]=([0, 1], \oplus, *, 0)$ is an MV-algebra and $x\odot y=\max \{0, x+y-1\}$.
	
	The rational numbers in $[0, 1]$
	and for each $n\geq 2$, the $n$-element set
	$$L_{n}=\{0, \frac{1}{n-1}, \frac{2}{n-1}, \cdots, \frac{n-2}{n-1}, 1\},$$
	yield examples of subalgebras of $[0, 1]$.
\end{example}

The MV-algebra $[0, 1]$ is important because it generates
the variety of all MV-algebras, and Chang Completeness Theorem  says that an equation holds in $[0, 1]$ if and only if it holds in every MV-algebra. Also, every $n$-element MV-chain is isomorphic to $L_{n}$ \cite{mu}.

For any $x, y \in A$, write 
$x\leq y$ for $ x^{*} \oplus y=1$.
It is well-known that  $\leq$
is a partial order on
$A$, called \name{the natural order} of $A$. Moreover, the natural order determines a structure
of bounded distributive lattice on $A$, with $0$ and $1$
are respectively the bottom and the top element, and
$$x\vee y=(x\odot y^{*})\oplus y~~\textrm{and} ~~ x\wedge y=x\odot (x^{*}\oplus y).$$
Clearly $(x\vee y)^{*}=x^{*}\wedge y^{*}$ and
$(x\wedge y)^{*}=x^{*}\vee y^{*}$.
We denote this bounded distributive lattice $(A, \vee, \wedge, 0, 1)$ by $\mathrm{L}(A)$,
and write $\mathrm{B}(A)$ as the set of all complemented elements of $\mathrm{L}(A)$.
An MV-chain is a linearly ordered MV-algebra.

\begin{lemma}\label{lem:1}
\cite[Lemmas 1.1.2 and 1.1.4]{mu} ~Let  $x, y \in A$. Then the following conditions are
	equivalent:
	\begin{enumerate}
		\item    $x\leq y$;
		\item  $x^{*} \oplus y=1$;
		\item  $x\odot y^{*}=0$;
		\item $y=x\vee y$;
		\item there exists an element $z\in A$ such that $x\oplus z=y$;
		\item $y^{*}\leq x^{*}$;
		\item  $x\oplus z\leq y\oplus z$ for all $z\in A$,
		\item  $x\odot z\leq y\odot z$ for all $z\in A$.
	\end{enumerate}
\end{lemma}

\begin{lemma}\label{lem:2}
\emph{(\cite[Proposition 1.1.6 and Lemma 1.6.1]{mu} and \cite[Theorem 1.15]{cc})}~Let $x, y, z\in A$. Then the following statements are true:
	\begin{enumerate}
		\item   $x\odot y\leq x\wedge y\leq x, y\leq x\vee y\leq x\oplus y$;\label{it:2.41}
		\item    $x\oplus (y\wedge z)=(x\oplus y)\wedge (x\oplus z)$;\label{it:2.42}
		\item   $x\odot (y\vee z)=(x\odot y)\vee (x\odot z)$;\label{it:2.43}
		\item $x\oplus y=y$ iff $x\odot y=x$ iff $x\wedge y^{*}=0$ iff $x^{*}\vee y=1$;\label{it:2.44}
		\item  If $x\odot y =x\odot z$ and $x\oplus y=x\oplus z$, then $y=z$.\label{it:2.45}
	\end{enumerate}
\end{lemma}

\begin{example}\label{exa:2}
	For any Boolean algebra  $(A, \vee, \wedge, -, 0, 1)$, the structure
	$(A, \vee, -, 0)$ is an MV-algebra, where $\vee, -$ and $0$ denote, respectively, the join,
	the complement and the smallest element in $A$.
\end{example}

Boolean algebras form a subvariety of the variety of MV-algebras. They
are precisely the MV-algebras satisfying the additional equation $x \oplus x = x$.
An MV-algebra $A$ is a Boolean algebra if and only if $A=\mathrm{B}(A)$.

\begin{lemma}\label{lem:3}
\cite[Theorem 1.5.3]{mu} For any $x\in A$, the following conditions are equivalent:
	\begin{enumerate}
		\item    $x\in \mathrm{B}(A)$;
		\item   $x\vee x^{*}=1$;
		\item  $x\wedge x^{*}=0$;
		\item   $x\oplus x=x$;
		\item   $x\odot x=x$;
		\item   $x^{\ast}\oplus x^{*}=x^{*}$;
		\item   $x^{*}\odot x^{*}=x^{*}$;
		\item    $x^{*}\in \mathrm{B}(A)$;
		\item   $x\oplus y=x\vee y$ for all $y \in A$;
		\item   $x\odot y=x\wedge y$ for all $y \in A$.
	\end{enumerate}
\end{lemma}

\subsection{Some notions of graphs}\label{subsec:2.2}	

A \name{graph} $\Gamma$ is an ordered pair $(V(\Gamma), E(\Gamma))$ consisting of a set $V (\Gamma)$ of vertices and a set $E(\Gamma)$, disjoint from $V(\Gamma)$, of edges, together with an incidence function $\psi_{\Gamma}$ that associates with each edge of $\Gamma$ an unordered pair of (not necessarily distinct) vertices of $\Gamma$. If $e$ is an edge and $u$ and $v$ are vertices such that $\psi_{\Gamma}(e)=\{u, v\}$, then $e$ is said to join $u$ and $v$, and the vertices $u$ and $v$ are called the ends of $e$. The ends of an edge are said to be incident with the
edge, and vice versa.

Two vertices which are incident with a common edge are
\name{adjacent}, 
an edge with identical ends is called a \name{loop}, and an edge with distinct ends a
\name{link}. Two or more links with the same pair of ends are said to be \name{parallel edges}.
A graph is \name{simple} if it has no loops or parallel edges.
If two vertices $u$ and $v$   are  adjacent in a simple graph, then we use the usual notation 
$u ~-~ v$.


A graph $G$ is called a \name{subgraph} of a graph $\Gamma$ if $V(G)\subseteq V(\Gamma)$,
$E(G)\subseteq E(\Gamma)$, and $\psi_{G}$ is the restriction of $\psi_{\Gamma}$ to $E(G)$.
Suppose that $X$ is a nonempty subset of $V(\Gamma)$. The subgraph of $\Gamma$ whose vertex set is $X$ and whose edge set is the set of those edges of $\Gamma$ that have both ends in $X$ is called \name{the subgraph of $\Gamma$ induced by $X$} and is denoted by $\Gamma[X]$.

By starting with a disjoint union of two graphs
$G$ and $H$ and adding edges joining every vertex of $G$ to every vertex of $H$, one
obtains the \name{join} of $G$ and $H$, which is denoted by $G\vee H$.

In what follows we always assume  that graphs are  simple.
Some other notions about graph theory are list in the following:
\begin{enumerate}
	\item [$\bullet$]
	A graph is \name{connected} if, for every partition of its vertex set into two nonempty
	sets $X$ and $Y$, there is an edge with one end in $X$ and one end in $Y$; otherwise the
	graph is \name{disconnected}. In other words, a graph is disconnected if its vertex set can
	be partitioned into two nonempty subsets $X$ and $Y$ so that no edge has one end
	in $X$ and one end in $Y$.
	\item [$\bullet$]  A \name{complete graph}
	is a simple graph in which any two distinct vertices are adjacent, and we denote by $K_{n}$
	the $n$-vertex complete graph.
	\item [$\bullet$]  An \name{empty graph}  is a graph whose edge set is empty, we denote
	by $\emptyset_{n}$  the empty graph with $n$-vertex. Notice that $K_{1}=\emptyset_{1}$.
	\item [$\bullet$]  The graph with no vertices (and hence no edges) is the \name{null graph}.
	\item [$\bullet$] A \name{path} is a simple graph whose vertices can be arranged in a linear sequence in such a way that two vertices are adjacent if they are consecutive in the sequence, and are nonadjacent otherwise. The \name{length} of a path is the number of its edges.
	\item [$\bullet$]   The \name{distance} $\D_{\Gamma}(a, b)$ between a pair of vertices $a$ and $b$
	in a graph $\Gamma$ is the length of the shortest path between them.
	\item [$\bullet$]The \name{diameter} $\diam(\Gamma)$ of a graph $\Gamma$ is
	defined to be the supremum of the distances between any pair of vertices.
	\item [$\bullet$]
	A graph is \name{bipartite} if its vertex set can be partitioned into two subsets $X$ and $Y$ so that every edge has one end in $X$ and one end in $Y$; such a partition $(X, Y )$ is called
	a \name{bipartition} of the graph. We denote a bipartite graph
	$\Gamma$ with bipartition $(X, Y )$ by $\Gamma[X, Y ]$. If $\Gamma[X, Y ]$ is simple and every vertex in $X$ is joined to every vertex in $Y$, then $\Gamma$ is called \name{a complete bipartite graph}. When $|X|=m$ and $|Y|=n$, we  denote the complete bipartite graph
	$\Gamma[X, Y ]$
	with bipartition $(X, Y )$ by $K_{m, n}$. It is easy to see that
	$K_{m, n}\cong \emptyset_{m}\vee \emptyset_{n}$.
	\item[$\bullet$] A \name{cycle} on three or more vertices is a simple graph whose vertices can be arranged in a cyclic sequence, vertices are adjacent if they are consecutive in the sequence, and are nonadjacent otherwise. The \name{length} of a cycle is the number of its edges.
	\item[$\bullet$] The \name{girth} of a  graph $\Gamma$, denoted by $\gr(\Gamma)$, is the length of a shortest cycle in $\Gamma$. If $\Gamma$ contains no cycles, then $gr(\Gamma)=\infty$.
\end{enumerate}

For more notions about MV-algebras and graphs, one can refer, respectively,  to \cite{mu} and \cite{bon}.

\section{Ideals of MV-algebras}\label{sec:3}

\begin{definition}\label{def:2111}
	\cite[Definition 4.1]{cc}
	Let $I$ be a subset of $A$.
	$I$  is called \name{an ideal} of $A$ if it satisfies: (i) $0\in I$,
	(ii) if $a, b\in I$, then
	$a\oplus b\in I$, and (iii) if $b\in I$ and $a\leq b$, then $a\in I$.
\end{definition}

It is clear that $\{0\}$ and $A$ itself are ideals of $A$. An ideal $I$ is said to be \name{proper} if $I\neq A$. It is easy to see that an ideal $I$ is proper if
and only if $1\not\in I$.
We denote by $\Ideal(A)$ the set of all ideals of $A$. 
 An MV-algebra $A$ is \name{simple} if $A$ is nontrivial and  $\{0\}$ is its only proper ideal.



\begin{lemma}\label{lem:44}
\cite[Theorem 3.5.1]{mu}  For every MV-algebra $A$ the following conditions are equivalent:
\begin{enumerate}
 \item    $A$ is simple;
\item   $A$ is nontrivial and for every nonzero element $a\in A$ there is a positive integer $n$ such that $1=a\oplus a \oplus \cdots \oplus a ~(n ~~times)$
\item $A$ is isomorphic to a subalgebra of $[0, 1]$.
\end{enumerate}
\end{lemma}

\begin{example}\label{exa:20}
	\emph{\cite{cc}}~ Define the following families of formal symbols:
	$$\mathcal{C}_{0} = \{0, c, 2c, 3c, \cdots\}, \quad \mathcal{C}_{1} = \{1, c^{*}, (2c)^{*}, (3c)^{*}, \cdots\}$$
	where $0=0c$, $1=0^{*}$, $(kc)^{*}=1-kc$, and $(kc)^{**}=((kc)^{*})^{*}=kc$. Letting $+ ~(resp. -)$ be the ordinary sum (resp. minus) between integers,
	define the following binary operation $\oplus$ on
	$\mathcal{C}=\mathcal{C}_{0}\cup \mathcal{C}_{1}$:
	\begin{enumerate}
		\item [$\bullet$] ~   $nc \oplus mc=(n+m)c$
		\item [$\bullet$] ~ $ (nc)^{*} \oplus (mc)^{*}=1$
		\item [$\bullet$]  ~ $ nc \oplus (mc)^{*}=(mc)^{*} \oplus nc=
		\begin{cases}
		1,  & \textrm{if}~ m\leq n; \\
		((m-n)c)^{*},  & \textrm{if}~ m>n.
		\end{cases}$
	\end{enumerate}
	Then $(\mathcal{C}, \oplus, *)$ is an MV-chain with $0<c<2c<3c<\cdots < (n-1)c<nc< \cdots < (nc)^{*}< ((n-1)c)^{*}< \cdots <(3c)^{*}< (2c)^{*}< c^{*}< 1$,
	and $\mathcal{C}_{0}$ is an ideal of $\mathcal{C}$.
\end{example}

The distance function $d: A \times A \rightarrow A$ is defined by
$$ d(x, y) = (x\odot y^{*}) \oplus (y\odot x^{*}).$$
For all $x, y, z\in A$, the following well-known
properties hold \cite[Proposition 1.2.5]{mu}:
\begin{enumerate}
	\item [$\bullet$]    $d(x, x)=0, ~d(x, y)=d(y, x), ~d(x, y)=d(x^{*}, y^{*}),
	~d(x, 0)= x$, and $d(x, 1) =x^{*}$;
	\item [$\bullet$]  if $d(x, y) =0$, then $x=y$;
	\item [$\bullet$]  $d(x, z)\leq d(x, y)\oplus d(y, z)$;
	\item [$\bullet$] $d(x\oplus u, y\oplus v)\leq d(x, y)\oplus d(u, v)$;
\end{enumerate}

\begin{proposition}\label{pro:20}
\cite[Proposition 1.2.6]{mu}~
	Let $I\in \Ideal (A)$. Define the binary
	relation $\equiv_{I}$ on $A$ defined by $x\equiv_{I}y$ if and only if
	$d(x, y)\in I$. Then $\equiv_{I}$ is a congruence relation. Moreover, $I = \{x \in A |~ x \equiv_{I} 0\}$.
	
	Conversely, if $\equiv$ is a congruence on $A$, then $\{x\in A|~ x\equiv 0\}\in \Ideal (A)$
	and $x \equiv y$ if and only if $d(x, y) \equiv 0$. Therefore, the correspondence
	$I\mapsto \equiv_{I}$ is a bijection
	from  $\Ideal (A)$ onto the set of congruences on $A$.
\end{proposition}

Let $I\in \Ideal (A)$ and $x \in A$. The equivalence class of $x$ with respect to $\equiv_{I}$ will be denoted
by $x/I$ and the quotient set $A/\equiv_{I}$ by $A/I$. Since $\equiv_{I}$ is a congruence, defining on the set $A/I$ the operations
\begin{center}
	$(x/I)^{\divideontimes} = x^{*}/I$ \quad and \quad $x/I \oplus y/I = (x \oplus y)/I$,
\end{center}
the system $(A/I, \oplus, \divideontimes, 0/I)$ becomes an MV-algebra, called \name{the quotient algebra}
of $A$ by the ideal $I$. Moreover, the correspondence $x \mapsto x/I$ defines a homomorphism
$h_{I}$ from $A$ onto the quotient algebra $A/I$, which is called \name{the natural
homomorphism} from $A$ onto $A/I$ \cite{mu}.

The following lemma is foundational, and will be used often in the sequel.

\begin{lemma}\label{lem:9}
	Let $I\in \Ideal (A)$, and $x, y\in A$. Then $x/I=y/I$ if and only if
	 $x\odot y^{*}\in I$ and $y\odot x^{*}\in I$.
\end{lemma}
\begin{proof}
	Since $x\odot y^{*}, y\odot x^{*}\leq(x\odot y^{*}) \oplus (y\odot x^{*})=d(x, y)$,
	it follows by  Definition \ref{def:2111} and Proposition \ref{pro:20} that the lemma is true.
\end{proof}

For notational convenience, we list some notations in the following:
\begin{enumerate}
	\item[$\bullet$]   $|X|$ denotes the cardinality of a set $X$.
	\item[$\bullet$]   $X\backslash Y=\{x\in X ~|~x\not\in Y\}$ for any sets $X$ and $Y$.
	\item[$\bullet$]    For a subset $B$ of  $A$, we denote the set
	$\{b^{*}~| ~b\in B\}$ by $B^{*}$.  Then
	$B^{**}=(B^{*})^{*}=B$. 
\end{enumerate}

\begin{lemma}\label{lem:10}
	Let  $I$ be a proper ideal of  $A$. Then for any $a\in A$, the following statements are true:
	\begin{enumerate}
		\item  $a^{*}/I=(a/I)^{*}$. In particular,
		$0/I=I$ and $1/I=I^{*}$.
		\label{it:2.111}
		\item $|a/I|=|a^{*}/I|$.
		\label{it:2.112}
		\item  $|a/I|\geq 2$ if $I\neq \{0\}$.
			\label{it:2.113}
		\item  $a/I$ is closed under $\wedge$ and $\vee$.	\label{it:2.114}
	\end{enumerate}
\end{lemma}
\begin{proof}
	Assume that $I$ is a proper ideal of  $A$.  Let $a\in A$.

\mnoindent	
\mref{it:2.111}  For any $x\in a/I$, we have $d(x^{*}, a^{*})=d(x, a)\in I$, and so  $x^{*}\in a^{*}/I$. Hence
	$(a/I)^{*}\subseteq a^{*}/I$.
Also, for any $y\in a^{*}/I$, we have $d(y^{*}, a)=d(y, a^{*})\in I$, and so
	$y^{*}\in a/I$, i.e. $y\in (a/I)^{*}$. Thus
	$a^{*}/I\subseteq (a/I)^{*}$, and hence $a^{*}/I=(a/I)^{*}$.
	By Proposition \ref{pro:20}, we have $0/I=I$, and consequently $1/I=0^{*}/I=(0/I)^{*}=I^{*}$.
	
\mnoindent
\mref{it:2.112}  follows from the fact that $\psi: a/I\rightarrow a^{*}/I; x\mapsto x^{*}$
	is a bijective map.
	
\mnoindent
\mref{it:2.113} Assume that $I\neq \{0\}$.
	Then $|1/I|=|0/I|=|I|\geq 2$ by \ref{it:2.111} and \ref{it:2.112}.  If there exists $a/I\not\in \{0/I, 1/I\}$ such that
	$|a/I|=1$, then $a/I=\{a\}$, and so $a^{*}/I=\{a^{*}\}$ by \ref{it:2.112}.
		Since $I\neq \{0\}$, there exists
	$u\in I$ such that $u\neq 0$. We have $(u\oplus a)/I=u/I\oplus a/I=0/I\oplus a/I=a/I$ and $$(u\oplus a^{*})/I=u/I\oplus a^{*}/I=0/I\oplus a^{*}/I=a^{*}/I,$$ which implies that $u\oplus a\in a/I=\{a\}$ and  $u\oplus a^{*}\in a^{*}/I=\{a^{*}\}$, so
	 $u\oplus a=a$ and  $u\oplus a^{*}=a^{*}$. It follows by Lemma \ref{lem:2} \ref{it:2.41} and \ref{it:2.44}
	that $u\leq u\oplus a=a$
	and $u\wedge a=0$, but then we get that
	$u=u\wedge a =0$, a contradiction. Therefore $|a/I|\geq 2$, and \ref{it:2.113} holds.
	
\mnoindent
\mref{it:2.114} Let $b, c\in a/I$. Then  $b\odot c^{*}\in I$ and $c\odot b^{*}\in I$ by Lemma \ref{lem:9}.
	Also, we have
	$(b\wedge c)\odot b^{*}\leq b\odot b^{*}=0$
	and $b\odot (b\vee c)^{*}=b\odot(b^{*}\wedge c^{*})\leq b^{*}\odot b=0$  by Lemma \ref{lem:1}; and
	 we get
	$$b\odot(b\wedge c)^{*} =b\odot (b^{*}\vee c^{*})=(b\odot b^{*})\vee (b\odot c^{*})=0\vee (b\odot c^{*})=b\odot c^{*}$$
	and $(b\vee c)\odot b^{*}=(b\odot b^{*})\vee (c\odot b^{*})= c\odot b^{*}$ by Lemma \ref{lem:2}.
	Thus we obtain that
	$$d(b\wedge c, b)=((b\wedge c)\odot b^{*})\oplus ( b\odot(b\wedge c)^{*})=b\odot c^{*}\in I,$$
	and $d(b\vee c, b)=((b\vee c)\odot b^{*})\oplus ( b\odot(b\vee c)^{*})=c\odot b^{*}\in I,$
	which implies that $b\wedge c, b\vee c \in b/I=a/I$, and hence  \ref{it:2.114} holds.
\end{proof}

\begin{proposition}\label{pro:15}
	Let  $I$ be a proper ideal of $A$. If $I$ is finite, then $|a/I|=|I|$ for any $a\in A$.
	In particular, if $A$ is finite, then $|a/I|=|I|$ for any $a\in A$.
\end{proposition}
\begin{proof}
Let $I$ be a proper and finite ideal of  $A$.
	If $I= \{0\}$, then it is clear that $|a/I|=|I|=1$ for any $a\in A$. So we sssume that  $I\neq \{0\}$ and  
	$I= \{0, i_{1}, i_{2}, \cdots, i_{k}\} ~(k\geq 1)$. 
	Put $u=i_{1}\oplus i_{2}\oplus\cdots \oplus i_{k}$.
	Then $u\in I$ by Definition \ref{def:2111}, and $x\leq u$ for any $x\in I$.
	
	For any  $a\in A$, if $a/I\in \{0/I, 1/I\}$, then $|a/I|=|I|$ by Lemma \ref{lem:10}.
	If $a/I\not\in \{0/I, 1/I\}$,
	define 
	$\varphi: a/I\rightarrow I$ by $\varphi(x) = x\odot u$ for any $x\in a/I$, then 	$\varphi$ is a map, since $x\odot u\leq u$ and $u\in I$ imply that $x\odot u\in I$. Also
	$\varphi$ is injective. In fact,
	suppose $x, y\in a/I$ such that $\varphi(x)=\varphi(y)$. Then we have
	$x\odot u=y\odot u$, and
	$x\odot y^{*}, x^{*}\odot y\in I$ by Lemma \ref{lem:9}. Thus
	\begin{eqnarray}\label{eqn:211}
	(x^{*}\odot y)\odot u=x^{*}\odot (y\odot u)=x^{*}\odot (x\odot u)
	=(x^{*}\odot x)\odot u=0\odot u.
	\end{eqnarray}
	Notice that $u$ is the largest element in $I$ and $x^{*}\odot y\in I$. We have
	$u\leq (x^{*}\odot y)\oplus  u\leq u$, so
	\begin{eqnarray}\label{eqn:212}
	(x^{*}\odot y)\oplus  u=u=0\oplus u.
	\end{eqnarray}
	Both $(\ref{eqn:211})$ and $(\ref{eqn:212})$ imply that $x^{*}\odot y=0$ by Lemma \ref{lem:2} \ref{it:2.45}, and so $y\leq x$ by Lemma \ref{lem:1}.
	Similarly, we can get $x\leq y$. Thus $x=y$, and so $\varphi$ is  injective.
	Hence $|a/I|\leq |I|$, and consequently $a/I$ is a finite subset of $A$, since $I$ is finite.
	
	Let $a/I=\{a_{1}, a_{2}, \cdots, a_{m}\}$ and put
	$v=a_{1}\wedge a_{2}\wedge \cdots \wedge a_{m}$. We have $v\in a/I$ by Lemma \ref{lem:10} \ref{it:2.114}, and so
	$v$ is the least element in $a/I$.
	Define 
	$\delta: I\rightarrow a/I$ by $\delta(x) = x\oplus v$ for any $x\in I$. Then $\delta$ is a map. In fact,
	for any $x\in I$, we have $$(x\oplus v)/I=x/I\oplus v/I=0/I\oplus v/I=v/I=a/I,$$
	and so $x\oplus v\in a/I$. Hence $\delta$ is a map.
	Also, $\delta$ is injective. Indeed, suppose $x, y\in I$ such that $\delta(x)=\delta(y)$, i.e.,
	$x\oplus v=y\oplus v$. Then
	\begin{eqnarray}\label{eqn:213}
	(x^{*}\oplus y)\oplus v=x^{*}\oplus (y\oplus v)=x^{*}\oplus (x\oplus v)=(x^{*}\oplus x)\oplus v
	=1\oplus v.
	\end{eqnarray}
	Since $x, y\in I$, we have $x/I=y/I=0/I$ and $x^{*}/I=(x/I)^{\divideontimes}=1/I$, and so
	$$((x^{*}\oplus y)\odot v)/I=
	(x^{*}/I\oplus y/I)\odot v/I=(1/I\oplus 0/I)\odot v/I=v/I=a/I.$$
	It follows that $(x^{*}\oplus y)\odot v \in a/I$.
	Notice that $v$ is the least element in $a/I$, we obtain that $v\leq (x^{*}\oplus y)\odot v\leq v$, and thus
	\begin{eqnarray}\label{eqn:214}
	(x^{*}\oplus y)\odot v=v=1\odot v.
	\end{eqnarray}
	Both $(\ref{eqn:213})$ and $(\ref{eqn:214})$ imply that $x^{*}\oplus y=1$ by Lemma \ref{lem:2} \ref{it:2.45},
	and so  $x\leq y$ by Lemma \ref{lem:1}.
	Similarly, we can get that $y\leq x$, and hence $x=y$. Therefore, $\delta$ is injective and
	consequently $|I|\leq |a/I|$. 
	
Summarizing the above arguments, we obtain
$|a/I|=|I|$.	
\end{proof}

Corollary \ref{cor:3.77}, closely analogous to Lagrange's Theorem for finite groups, follows immedialtely 
from Proposition \ref{pro:15}.
\begin{corollary}\label{cor:3.77}
	Let  $I\in \Ideal (A)$. Then $A$ is finite if and only if $I$ and $A/I$ are finite. Moreover, $|A|=|I|\cdot |A/I|$
	when  $A$ is finite.
\end{corollary}

Let $\Omega$ be an index set. The direct product  $\prod_{i\in \Omega}A_{i}$ of a family $\{A_{i}\}_{i\in \Omega}$ of MV-algebras
is the MV-algebra obtained by endowing the set-theoretical cartesian product of the $A_{i}$'s with the pointwise MV-operations.
The zero element $\mathbf{0}$ of $\prod_{i\in \Omega}A_{i}$ is the function $i\in \Omega\mapsto 0_{i}\in A_{i}$, and the element $\mathbf{1}$
of $\prod_{i\in \Omega}A_{i}$ is the function $i\in \Omega\mapsto 1_{i}\in A_{i}$ (see \cite{mu}).
Let $f, g\in \prod_{i\in \Omega}A_{i}$. It is clear that
$f\leq g$ in $\prod_{i\in \Omega}A_{i}$ if and only if
$f(i)\leq g(i)$ in $A_{i}$ for each $i\in \Omega$.

When $\Omega=\{1, 2, \cdots, n\}$, we denote the direct product of $\{A_{i}\}_{i\in \Omega}$ by
$\prod_{i=1}^{n}A_{i}$ or $A_{1}\times A_{2}\times \cdots \times A_{n}$.
We know that $L_{2}\times L_{2}\cong B_{4}$, where $B_{4}$ is the $4$-element Boolean algebra (see \cite{gan}).

\begin{theorem}\label{the:1000}
	Let $\Omega$ be an index set, and $A=\prod_{i\in \Omega}A_{i}$
	be the direct product  of a family $\{A_{i}\}_{i\in \Omega}$ of MV-algebras.
	If $I_{i}\in \Ideal(A_{i})$ for each $i\in \Omega$, then the cartesian  product $\prod_{i\in \Omega}I_{i}$ of sets $I_{i}$  is an ideal of $A$, and $A/(\prod_{i\in \Omega}I_{i})\cong \prod_{i\in \Omega}(A_{i}/I_{i})$.
\end{theorem}

\begin{proof}
	Let $\Omega, A, A_{i}$ and $I_{i}$ be as given.
	It is routine to verify that $\prod_{i\in \Omega}I_{i}$ is an  ideal of $A$. 
	
	Define $\psi: A/(\prod_{i\in \Omega}I_{i})\rightarrow \prod_{i\in \Omega}(A_{i}/I_{i})$ by
	$\psi(f/(\prod_{i\in \Omega}I_{i}))= \prod_{i\in \Omega} (f(i)/I_{i})$ for any
	$f/(\prod_{i\in \Omega}I_{i})\in A/(\prod_{i\in \Omega}I_{i})$.
	Let $f, g\in A$.
	We have
	\begin{eqnarray*}
		&&               f/(\prod_{i\in \Omega}I_{i})=g/(\prod_{i\in \Omega}I_{i})\\
		&\Leftrightarrow& f\odot g^{*}\in \prod_{i\in \Omega}I_{i} ~ \textrm{and}~
		g\odot f^{*}\in \prod_{i\in \Omega}I_{i}
		\quad ~\mbox{(by Lemma \ref{lem:9})}\\
		&\Leftrightarrow& (\forall i\in \Omega) ~f(i)\odot g(i)^{*}\in I_{i}  ~ \textrm{and} ~
		g(i)\odot f(i)^{*}\in I_{i} \\
		&\Leftrightarrow& (\forall i\in \Omega) ~ f(i)/I_{i}=g(i)/I_{i} \quad~\mbox{(by Lemma \ref{lem:9})}\\
		&\Leftrightarrow& \prod_{i\in \Omega}( f(i)/I_{i})=\prod_{i\in \Omega}( g(i)/I_{i})\\
		&\Leftrightarrow& \psi(f/(\prod_{i\in \Omega}I_{i}))=\psi(g/(\prod_{i\in \Omega}I_{i})).
	\end{eqnarray*}
	and
	$$\psi( f/(\prod_{i\in \Omega}I_{i}) \oplus g/(\prod_{i\in \Omega}I_{i}))
	=\psi((f\oplus g)/(\prod_{i\in \Omega}I_{i}))$$$$
	= \prod_{i\in \Omega} (f\oplus g)(i)/I_{i}
	= \prod_{i\in \Omega} (f(i)\oplus g(i))/I_{i}
	= \prod_{i\in \Omega} (f(i)/I_{i}\oplus g(i)/I_{i})$$$$
	= \prod_{i\in \Omega} (f(i)/I_{i})\oplus \prod_{i\in \Omega} (g(i)/I_{i})
	= \psi(f/(\prod_{i\in \Omega}I_{i}))\oplus\psi(g/(\prod_{i\in \Omega}I_{i}))$$
	and
	$$\psi((f/(\prod_{i\in \Omega}I_{i}))^{\divideontimes})=\psi(f^{*}/(\prod_{i\in \Omega}I_{i}))
	=\prod_{i\in \Omega} (f(i)^{*}/I_{i})$$
	$$=\prod_{i\in \Omega} (f(i)/I_{i})^{\divideontimes}=(\prod_{i\in \Omega} f(i)/I_{i})^{\divideontimes}=
	\psi((f/(\prod_{i\in \Omega}I_{i})))^{\divideontimes}.$$
	Hence $\psi$ is an injective homomorphism.
	Also, it is obvious that $\psi$ is surjective. Thus $\psi$ is  an isomorphism,
	and hence  $A/(\prod_{i\in \Omega}I_{i})\cong \prod_{i\in \Omega}(A_{i}/I_{i})$.
\end{proof}

\begin{remark}
	Let $\Omega$ be an index set, and $A=\prod_{i\in \Omega}A_{i}$
be the direct product  of a family $\{A_{i}\}_{i\in \Omega}$ of MV-algebras. Theorem \ref{the:1000} says that $\{\prod_{i\in \Omega}I_{i} ~| ~(\forall i\in \Omega)~ I_{i}\in \Ideal (A_{_{i}}) \}\subseteq \Ideal(A)$. But generally,
$\{\prod_{i\in \Omega}I_{i} ~|~(\forall i\in \Omega)~ I_{i}\in \Ideal (I_{_{i}})\}\neq \Ideal(A)$. For example, if $\Omega$ is infinite, let 
\begin{center}	
$P=\{f\in A~|~\Supp (f) ~is ~finite \}$, where $\Supp(f)=\{i\in \Omega ~|~ f(i)\neq 0\}$.
\end{center}
 Then it is easy to see that $P$ is a proper ideal of $A$, but $P\not\in \{\prod_{i\in \Omega}I_{i} ~|~(\forall i\in \Omega)~ I_{i}\in \Ideal (I_{_{i}})\}.$
Indeed, if $P=\prod_{i\in \Omega}I_{i}$, where
$I_{i}\in \Ideal(A_{_{i}})$, then for each $i\in \Omega$ and any $z\in A_{i}$, let
 $$\displaystyle g(j)= \left\{ \begin{array}{ll}
0, & \textrm{if} \quad j\in \Omega\backslash \{i\} \\
z, & \textrm{if} \quad j= i
\end{array} \right.$$
Then $g\in P$, and so $z=g(i)\in I_{i}$, which implies that $I_{i}=A_{i}$, and hence $P=A$, a contradiction.
\end{remark}

\begin{proposition}\label{cor:10000}
	Let $\Omega$ be a finite index set, and $A=\prod_{i\in \Omega}A_{i}$
	be the direct product  of a family $\{A_{i}\}_{i\in \Omega}$ of MV-algebras.
Then  $\{\prod_{i\in \Omega}I_{i} ~| ~(\forall i\in \Omega)~ I_{i}\in \Ideal (A_{_{i}}) \}= \Ideal(A)$.
\end{proposition}

\begin{proof}	Let $\Omega, A$ and $A_{i}$ be as given. Then  $\{\prod_{i\in \Omega}I_{i} ~| ~(\forall i\in \Omega)~ I_{i}\in \Ideal (A_{_{i}}) \}\subseteq \Ideal(A)$
	by Theorem \ref{the:1000}.
	Conversely, let $K\in \Ideal(A)$. For each $i\in \Omega$,
	put
	$$I_{i}=\{a_{i} \in A_{i}~|~(\exists~ f\in K)~ a_{i}=f(i)\}.$$
	Then $I_{i}$ is an ideal of $A_{i}$. Indeed, we have  $0\in I$, since $\mathbf{0}\in K$. Also,
	for any $a_{i}, b_{i}\in I_{i}$ and any $c\in A_{i}$ with $c\leq a_{i}$, there exist
	$f, g\in K$ such that $a_{i}=f(i)$ and $b_{i}=g(i)$, so
	$a_{i}\oplus b_{i}=f(i)\oplus g(i)=(f\oplus g)(i)$. It follows that $a_{i}\oplus b_{i}\in I_{i}$,
	since $f\oplus g\in K$. To prove that $c\in I_{i}$
	Let $$\displaystyle h(j)= \left\{ \begin{array}{ll}
	f(j), & \textrm{if} \quad j\in \Omega\backslash \{i\} \\
	c, & \textrm{if} \quad j= i
	\end{array} \right.$$
	Then $h\in \prod_{i\in \Omega}A_{i}=A$ and $h\leq f$, since $h(i)=c\leq a_{i}=f(i)$ and $h(j)=f(j)$ for any $j\in \Omega\backslash \{i\}$.
	It follows that $h\in K$, and so $c\in I_{i}$. Hence
	$I_{i}$ is an ideal of $A_{i}$.
	
	Finally, we claim that $K=\prod_{i\in \Omega} I_{i}$. In fact,  it is clear that
	$K\subseteq \prod_{i\in \Omega} I_{i}$. To prove the opposite inclusion, let $w\in \prod_{i\in \Omega} I_{i}$.
Then for each $i\in \Omega$, we have $w(i)\in I_{i}$, and so there exists
	$f^{i}\in K$  such that $f^{i}(i)=w(i)$. 
	Let $f=\oplus_{i\in \Omega}f^{i}$. Notice that $\Omega$ is finite and $w(i)\leq f(i)$, we have $f$ is well-defined, $f\in K$ and 
	$w\leq f$, which implies that $w\in K$, since $K$ is an ideal of $A$. Hence $\prod_{i\in \Omega} I_{i}\subseteq K$, and therefore $K=\prod_{i\in \Omega} I_{i}$. Consequently,  $\{\prod_{i\in \Omega}I_{i} ~| ~(\forall i\in \Omega)~ I_{i}\in \Ideal (A_{_{i}}) \}= \Ideal(A)$ as required.
\end{proof}

From Theorem \ref{the:1000}, we can immediately get
\begin{corollary}\label{cor:6}
	Let $n$ be a positive integer,  $\Omega=\{{1, 2, \cdots, n}\}$ and $A=\prod_{i=1}^{n}A_{i}$ be the direct product of  MV-algebras
	$A_{i}~ (i\in \Omega)$. Then $A/I_{i}\cong \prod_{j\neq i}A_{j}$ where
	$I_{i}=\{0\}\times\cdots\times \{0\}\times A_{i}\times\{0\}\times\cdots\times \{0\}$ for each $i\in \Omega$. 
\end{corollary}

\begin{example}\label{exa:10}
	Let $M=\{0, a, b, c, d, e, f, g, h, i, j, 1\}$ and operations $\oplus$,
	$*$ and $\odot$ be defined as follows
	$$
	\begin{tabular}{c|cccccccccccc}
	$*$         & $0$ & $a$  &$b$ & $c$ & $d$   &$e$ & $f$ & $g$  &$h$  & $i$ & $j$  &$1$ \\ \hline
	& $1$ & $j$  &$i$ & $h$ & $g$  &$f$  & $e$ & $d$  &$c$  & $b$ & $a$  &$0$ \\
	\end{tabular}
	$$
	\vspace{0.03cm}
	$$
	\begin{tabular}{c|cccccccccccc}
	$\oplus$     & $0$ & $a$  &$b$  & $c$ & $d$  &$e$ & $f$ & $g$  &$h$  & $i$ & $j$  &$1$\\ \hline
	$0$  & $0$ & $a$  &$b$  & $c$ & $d$  &$e$ & $f$ & $g$  &$h$  & $i$ & $j$  &$1$  \\
	$a$  & $a$ & $b$  & $b$ & $d$ & $e$  &$e$ & $g$ & $h$  &$h$  & $j$ & $1$  &$1$\\
	$b$  & $b$ & $b$  & $b$ & $e$ & $e$  &$e$ & $h$ & $h$  &$h$  & $1$ & $1$  &$1$\\
	$c$  & $c$ & $d$  &$e$  & $c$ & $d$  &$e$ & $i$ & $j$  &$1$  & $i$ & $j$  &$1$  \\
	$d$  & $d$ & $e$  & $e$ & $d$ & $e$  &$e$ & $j$ & $1$  &$1$  & $j$ & $1$  &$1$\\
	$e$  & $e$ & $e$  & $e$ & $e$ & $e$  &$e$ & $1$ & $1$  &$1$  & $1$ & $1$  &$1$\\
	$f$  & $f$ & $g$  &$h$  & $i$ & $j$  &$1$ & $f$ & $g$  &$h$  & $i$ & $j$  &$1$  \\
	$g$  & $g$ & $h$  & $h$ & $j$ & $1$  &$1$ & $g$ & $h$  &$h$  & $j$ & $1$  &$1$\\
	$h$  & $h$ & $h$  & $h$ & $1$ & $1$  &$1$ & $h$ & $h$  &$h$  & $1$ & $1$  &$1$\\
	$i$  & $i$ & $j$  &$1$  & $i$ & $j$  &$1$ & $i$ & $j$  &$1$  & $i$ & $j$  &$1$  \\
	$j$  & $j$ & $1$  & $1$ & $j$ & $1$  &$1$ & $j$ & $1$  &$1$  & $j$ & $1$  &$1$\\
	$1$  & $1$ & $1$  & $1$ & $1$ & $1$  &$1$ & $1$ & $1$  &$1$  & $1$ & $1$  &$1$\\
	\end{tabular}
	$$
		\vspace{0.03cm}
	$$
	\begin{tabular}{c|cccccccccccc}
	$\odot$      & $0$ & $a$  &$b$  & $c$ & $d$  &$e$ & $f$ & $g$  &$h$  & $i$ & $j$  &$1$\\ \hline
	$0$  & $0$ & $0$  &$0$  & $0$ & $0$  &$0$ & $0$ & $0$  &$0$  & $0$ & $0$  &$0$  \\
	$a$  & $0$ & $0$  & $a$ & $0$ & $0$  &$a$ & $0$ & $0$  &$a$  & $0$ & $0$  &$a$\\
	$b$  & $0$ & $a$  & $b$ & $0$ & $a$  &$b$ & $0$ & $a$  &$b$  & $0$ & $a$  &$b$\\
	$c$  & $0$ & $0$  &$0$  & $c$ & $c$  &$c$ & $0$ & $0$  &$0$  & $c$ & $c$  &$c$  \\
	$d$  & $0$ & $0$  & $a$ & $c$ & $c$  &$d$ & $0$ & $0$  &$a$  & $c$ & $c$  &$d$\\
	$e$  & $0$ & $a$  & $b$ & $c$ & $d$  &$e$ & $0$ & $a$  &$b$  & $c$ & $d$  &$e$\\
	$f$  & $0$ & $0$  &$0$  & $0$ & $0$  &$0$ & $f$ & $f$  &$f$  & $f$ & $f$  &$f$  \\
	$g$  & $0$ & $0$  & $a$ & $0$ & $0$  &$a$ & $f$ & $f$  &$g$  & $f$ & $f$  &$g$\\
	$h$  & $0$ & $a$  & $b$ & $0$ & $a$  &$b$ & $f$ & $g$  &$h$  & $f$ & $g$  &$h$\\
	$i$  & $0$ & $0$  &$0$  & $c$ & $c$  &$c$ & $f$ & $f$  &$f$  & $i$ & $i$  &$i$  \\
	$j$  & $0$ & $0$  & $a$ & $c$ & $c$  &$d$ & $f$ & $f$  &$g$  & $i$ & $i$  &$j$\\
	$1$  & $0$ & $a$  & $b$ & $c$ & $d$  &$e$ & $f$ & $g$  &$h$  & $i$ & $j$  &$1$\\
	\end{tabular}
	$$
	
	Since $\psi: L_{2}\times L_{2}\times L_{3}\rightarrow M$;
	$(0, 0, 0)\mapsto 0, (0, 0,\frac{1}{2})\mapsto a$,
	$(0, 0, 1)\mapsto b$,
	$(0, 1, 0)\mapsto c,
	(0, 1, \frac{1}{2})\mapsto d, (0, 1, 1)\mapsto e, (1, 0, 0)\mapsto f, (1, 0,\frac{1}{2})\mapsto g$, $(1, 0, 1)\mapsto h, (1, 1, 0)\mapsto i,
	(1, 1, \frac{1}{2})\mapsto j, (1, 1, 1)\mapsto 1$
	is an isomorphism,   $M$ is an MV-algebra and  $M\cong L_{2}\times L_{2}\times L_{3}$ by \cite[Example 2.10]{GY2}.
	
	Let $I=\{0, a, b\}$ and $J=\{0, c\}$. Then by Corollary \ref{cor:6} both $I$ and $J$ are ideals of $M$,
	$M/I=\{0/I, c/I, h/I, 1/I\}\cong L_{2}\times L_{2}\cong B_{4}$, and
	$M/J=\{0/J, a/J, b/J, i/J, j/J, 1/J\}\cong L_{2}\times L_{3},$
	where
	$0/I=I, c/I=\{c, d, e\}, h/I=\{f, g, h\}$, $1/I=\{i, j, 1\}$; and
	$0/J=J, a/J=\{a, d\}, b/J=\{b, e\}, i/J=\{i, f\}, j/J=\{j, g\}$, $1/J=\{h, 1\}$.
\end{example}

\section{The girth of the zero-divisor graph $\Gamma(A)$}\label{sec:4}

Recall that the zero-divisor graph
$\Gamma (A)$ of $A$ is the simple graph whose  vertex set $V(\Gamma (A))$  is the set of
all non-zero zero-divisors of $A$, such that two distinct such elements $a, b$ form an edge
precisely when $a\odot b = 0$. In
\cite{gan}, Gan and Yang discussed the diameter of $\Gamma (A)$, and use $\Gamma (A)$ to characterize $A$ whose cardinality up to $7$.
In this section, we will investigate the girth of $\Gamma (A)$, and use them to classify all MV-algebras. 

Remember that the girth of a  graph $\Gamma$, denoted by $\gr(\Gamma)$, is the length of a shortest cycle in $\Gamma$. If $\Gamma$ contains no cycles, then $gr(\Gamma)=\infty$.

\begin{lemma}\label{lem:5}
	\cite[Lemma 4]{gan}~
	Let  $\Fix_{*}(A)=\{a\in A ~|~  a^{*}=a\}$. Then
	$|\Fix_{*}(A)|\leq 1$.
\end{lemma}

\begin{proposition}\label{pro:4.20}
	Let $A$ be  an MV-algebra, and $A\cong A_{1}\times A_{2} $, where $A_{1}, A_{2}$ are nontrivial MV-algebras. 
	\begin{enumerate}
		\item  	If $|A_{1}|\geq 3$ and $|A_{2}|\geq 3$, then  	
		$ \gr(\Gamma(A))=3$.
		\label{it:4.201}
		\item   	If  $|A_{1}|=2$ and $|A_{2}|\geq 4$, then  	
		$ \gr(\Gamma(A))=3$.
		\label{it:4.202}	
		\item 	If $|A_{1}|=2$ and $|A_{2}|=2$, or $|A_{1}|=2$ and $|A_{2}|=3$, then  	
		$ \gr(\Gamma(A))=\infty$.
		\label{it:4.203}
	\end{enumerate}	
\end{proposition}
\begin{proof}
	Let $ A, A_{1}$ and $A_{2}$ be as given.
	
	If $|A_{1}|\geq 3$ and $|A_{2}|\geq 3$, 	
	fix $a\in A_{1}\backslash\{0, 1\}$ and $b\in A_{2}\backslash\{0, 1\}$, then $(a, 0)~-~(a^{*}, b) ~-~ (0, b^{*}) ~-~(a, 0)$ is a triangle in  $\Gamma(A_{1}\times A_{2} )$. Hence 
	$\gr(\Gamma (A))=\gr (\Gamma( A_{1}\times A_{2})=3$.
	
	If $|A_{1}|= 2$ and $|A_{2}|\geq 4$, then
	there exist	
	$b\in A_{2}\backslash\{0, 1\}$ such that $b\neq b^{*}$ by Lemma \ref{lem:5}, and so $(1, 0)~-~(0, b)~-~ (0, b^{*}) ~-~ (1, 0)$ is a triangle in  $\Gamma(A_{1}\times A_{2} )$. Hence 
	$\gr(\Gamma(A))=\gr (\Gamma( A_{1}\times A_{2})=3$.
	
	If  $|A_{1}|= 2$ and $|A_{2}|= 2$, then $A\cong B_{4}$, where $B_{4}$ is the $4$-element Boolean algebra, and $\Gamma(A)\cong K_{2}$ by \cite[Theorem 3]{gan}. Hence $\gr(\Gamma(A))=\infty$.
	
	If  $|A_{1}|= 2$ and $|A_{2}|= 3$, then $A\cong L_{2}\times L_{3}$, and $\Gamma(A)$ is a path of length $3$ by \cite[Examples 3, 4]{gan}. Hence 
	$ \gr(\Gamma(A))=\infty$.
\end{proof}

Recall that an MV-algebra $A$ is \name{directly indecomposable} if $A$ is not isomorphic to a direct product of two nontrivial MV-algebras.
An MV-algebra $A$ is directly indecomposable if and only if $\mathrm{B}(A)=\{0, 1\}$, i.e, $A$ has only two idempotent: $0$ and $1$ (see \cite{Cig} or \cite[Theorem 8]{Bel}).

\begin{proposition}\label{pro:4.21}
	Let $A$ be  a directly indecomposable MV-algebra. If $|A|\geq 6$, then 	$ \gr(\Gamma(A))=3$.
\end{proposition}
\begin{proof}
	Let $A$ be  a directly indecomposable MV-algebra with $|A|\geq 6$. 
If $|A|=6$, then $A\cong L_{6}$, and it is easy to see that 	$\gr(\Gamma(A))=\gr(\Gamma(L_{6}))=3$. Assume that $|A|\geq 7$. By Lemma \ref{lem:5},
 there exist 	
	$a, b\in A\backslash\{0, 1\}$ such that $a, b$, $ a^{*}$ and $ b^{*}$ are mutually distinct elements in $A$ 
	
	If $a$ and $a^{*}$ are incomparable (i.e, $a\not\leq a^{*}$ and $a^{*}\not\leq a$), then 
	$a\wedge a^{*}\neq 0$ by Lemma \ref{lem:3}, since $a$ is not an idempotent. Thus
	$a~-~ (a\wedge a^{*})~-~ a^{*}~-~a$ is a triangle  in  $\Gamma(A)$.  Similarly, if  $b$ and $b^{*}$ are incomparable, then 	there is a triangle in $\Gamma(A)$.	  
	
	So, without loss of generality, assume that $a< a^{*}$(i.e, $a\leq a^{*}$ and $a\neq a^{*}$) and $b< b^{*}$. Then $a\odot a=0$ and $b\odot b=0$ by Lemma \ref{lem:1}. If $a$ and $b$ are comparable, say $a<b$, then by Lemma \ref{lem:1} again, we have $a<b<b^{*}<a^{*}$,  and so $a~-~b~-~ b^{*}~-~a$ is a triangle  in  $\Gamma(A)$.  If $a$ and $b$ are incomparable and $a\wedge b\neq 0$, then $a~-~a^{*}~- ~(a\wedge b)~-~a$ is a triangle  in  $\Gamma(A)$, since $a\odot (a\wedge b)\leq a\odot a=0$ and $a^{*}\odot (a\wedge b)\leq a^{*}\odot a=0$ by Lemma \ref{lem:1}.  
	
	If $a$ and $b$ are incomparable and $a\wedge b= 0$, then $a\odot b=0$, since
	 $a\odot b\leq a\wedge b$ by Lemma \ref{lem:2} \ref{it:2.41}. It follows by Lemma \ref{lem:1} that $a\leq b^{*}$ and $b\leq a^{*}$, so $a\vee b\leq b^{*}$ and
	 $a\vee b\leq a^{*}$. Hence $a\vee b \leq a^{*}\wedge b^{*}$, and thus
	 $a~-~b~- ~ (a^{*}\wedge b^{*})~-~a$ is a triangle  in  $\Gamma(A)$, since $a\odot  (a^{*}\wedge b^{*})\leq a\odot a^{*}=0$ and  $b\odot  (a^{*}\wedge b^{*})\leq b\odot b^{*}=0$ by Lemma \ref{lem:1}.  	
	
	Summarizing the above arguments, we get
	$\gr(\Gamma(A))=3$.
\end{proof}

\begin{theorem}\label{the:4.30}
	Let $A$ be an MV-algebra with $|A|\geq 3$. Then
	$$\gr (\Gamma(A))=
	\begin{cases}
	\infty,  & \textrm{if}~ A\cong L_{3}, ~A\cong L_{4}, ~A\cong L_{5}, A\cong B_{4}~~ or~~ A\cong L_{2}\times L_{3}; \\
	3,  & \textrm{otherwise}.
	\end{cases}$$
\end{theorem}
\begin{proof}
	If $|A|=3$, then $\Gamma(A)\cong K_{1}$, and so
	$\gr(\Gamma(A))=\infty$.
	If $|A|=4$, then $\Gamma(A)\cong K_{2}$ by Theorem \cite[Theorm 3]{gan}, and so
	$\gr(\Gamma(A))=\infty$.
		If $|A|=5$, then $A\cong L_{5}$ by \cite[Theorem 5]{gan}, and so $\Gamma(A)$ is a path of length $2$. Hence 
	$\gr(\Gamma(A))=\infty$.
	
	If $|A|=6$, then $A$ is isomorphic to $ L_{6}$ or $L_{2}\times L_{3}$ by \cite[Theorem 7]{gan}.
	If $A\cong L_{6}$, then $\gr(\Gamma(A))=3$ by Proposition \ref{pro:4.21}. If $A\cong L_{2}\times L_{3}$, then $\gr(\Gamma(A))=\infty$ by Proposition \ref{pro:4.20}.
	
	If $|A|=7$, then $A\cong L_{7}$
	by \cite[Theorem 8]{gan}, so
	$\gr(\Gamma(A))=3$ by Proposition \ref{pro:4.21}.
	
	If $|A|\geq 8$, then $\gr(\Gamma(A))=3$ by Propositions \ref{pro:4.20} and  \ref{pro:4.21}.
\end{proof}	

\section{The ideal-based zero-divisor graph $\Gamma_{I}(A)$}\label{sec:5}

In this section, we will introduce and study the ideal-based zero-divisor graph $\Gamma_{I}(A)$ of an MV-algebra $A$ with respect to an ideal $I$ of $A$. The diameter and the girth of $\Gamma_{I}(A)$ will be discussed. 

\subsection{Some basic properties of $\Gamma_{I}(A)$}\label{subsec:4.1}

Let $I$ be a proper ideal of an MV-algebra $A$. We define the ideal-based zero-divisor graph
$\Gamma_{I}(A)$ of $A$ with respect to $I$ as follows:
\begin{enumerate}
	\item [$\blacktriangleright$]   The vertex set
	$V(\Gamma_{I}(A))=\{x\in A\backslash I ~|~ (\exists~ y\in A\backslash I) ~x\odot y\in I\}$.
	\item [$\blacktriangleright$]   There exists an edge  between distinct vertices $x$
	and $y$ if and only if $x\odot y\in I$.
\end{enumerate}
It is clear that if $I=\{0\}$, then $\Gamma_{I}(A)=\Gamma(A)$.


\begin{proposition}\label{pro:11}
	Let $I$ be a proper ideal of an MV-algebra $A$. Then
	$V(\Gamma_{I}(A))=A\backslash (I\cup I^{*})$. In particular, $V(\Gamma (A))=A\backslash\{0, 1\}$.
\end{proposition}
\begin{proof}
	Assume that $x\in V(\Gamma_{I}(A))$.
	Then $x\in A\backslash I$ and  $x\odot y\in I$ for some $y\in A\backslash I$.
	If $x\in I^{*}$, then $x^{*}\in I^{**}=I$, and so
	$x^{*}\vee y=x^{*}\oplus (x\odot y)\in I$. Since $I$ is an ideal of $A$ and $y\leq x^{*}\vee y$,
	we get $y\in I$, a contradiction. Thus $x\in A\backslash (I\cup I^{*})$, and hence
	$V(\Gamma_{I}(A))\subseteq A\backslash (I\cup I^{*})$.
	
	To prove the opposite inclusion, let $a\in A\backslash (I\cup I^{*})$. Then it is easy to see that
	$a^{*}\in A\backslash (I\cup I^{*})$ and
	$a\odot a^{*}=0\in I$, so $a\in  V(\Gamma_{I}(A))$. Thus
	$A\backslash(I\cup I^{*})\subseteq V(\Gamma_{I}(A))$.
\end{proof}

\begin{corollary}\label{cor:5}
	Let  $I$ be a proper ideal of an MV-algebra $A$, and $a\in A$. Then $a\in V(\Gamma_{I}(A))$ iff $a/I\in V(\Gamma(A/I))$ iff $a/I\subseteq V(\Gamma_{I}(A))$, where $\Gamma(A/I)$ is the zero-divisor graph of the MV-algebra $A/I$.
\end{corollary}
\begin{proof}
	Notice that $0/I=I$, $1/I=I^{*}$ and  $V(\Gamma(A/I))=(A/I) \backslash \{0/I, 1/I\}$.
	The result follows immediately by Proposition \ref{pro:11}.
\end{proof}

By  Proposition \ref{pro:11} and Corollary \ref{cor:5} we  obtain that $\Gamma_{I} (A)$ is null if and only if $\Gamma (A/I)$ is null, and if and only if
$A/I\cong L_{2}$.
Remember in mind that $a\in V(\Gamma_{I}(A))$ iff $a\not\in I\cup I^{*}$.

Compare the next result with \cite[Lemma 2.1]{at1}.

\begin{theorem}\label{the:1}	Let  $I$ be a proper ideal of an MV-algebra $A$ with $|A/I|\geq 3$. Then
	$\Gamma_{I} (A)$  is connected with $\diam(\Gamma_{I}(A))\leq 3$.
\end{theorem}

\begin{proof}
	Let $a, b\in V(\Gamma_{I} (A))$ such that $a\neq b$. Then $a^{*}, b^{*}\in V(\Gamma_{I} (A))$ and $a^{*}\neq b^{*}$.
	By Lemma \ref{lem:5}, without loss of generality, we can assume that $b^{*}\neq b$.
	We will show that $\D_{\Gamma_{I} (A)}(a, b)\leq 3$. Consider the following cases:
	
	If  $a\odot b\in I$, then there exists an edge $a~ -~b$ in $\Gamma_{I} (A)$.
	
	If $a\odot b\not\in I$ but $a\odot b^{*}\in I$, then there exists a path $a ~-~ b^{*}~-~b$ in $\Gamma_{I} (A)$.
	
	If  $a\odot b\not\in I, a\odot b^{*}\not\in I$ but $a^{*}\odot b\in I$, then $a\neq a^{*}$ and
	there exists a path $a ~- ~a^{*} ~- ~b$ in $\Gamma_{I} (A)$.
	
	If  $a\odot b\not\in I, a\odot b^{*}\not\in I, a^{*}\odot b\not\in I$ but $a^{*}\odot b^{*}\in I$, then $a\neq a^{*}$, and
	there exists a path $a ~-~ a^{*}~-~b^{*}~-~b$ in $\Gamma_{I} (A)$.
	
	If  $a\odot b\not\in I, a\odot b^{*}\not\in I, a^{*}\odot b\not\in I$ and $ a^{*}\odot b^{*}\not\in I$, then
	$a^{*}\odot b^{*}\not\in I^{*}$, for otherwise, $a^{*}\odot b^{*}\in I^{*}$ implies that
	$b\leq a\oplus b=(a^{*}\odot b^{*})^{*}\in I$, and so $b\in I$, a contradiction. Thus we get $a^{*}\odot b^{*}\in V(\Gamma_{I}(A))$ by Proposition \ref{pro:11}, and so
	there exists a path $a ~- ~( a^{*}\odot b^{*})~- ~b$ in $\Gamma_{I} (A)$. 
	
	Therefore  $\D_{\Gamma_{I} (A)}(a, b)\leq 3$, and hence $\Gamma_{I} (A)$  is connected and $\diam(\Gamma_{I}(A))\leq 3$.
\end{proof}

The following example  shows that there  exists  proper ideals $I$ and
$J$ of an MV-algebra $M$ such that $\diam(\Gamma_{I}(M))=2$ and $\diam(\Gamma_{J}(M))=3$.

\begin{example}\label{exm:1000}
	Let $M=\{0, a, b, c, d, e, f, g, h, i, j, 1\}$ be the MV-algebra in Example
	\ref{exa:10}, and  $I=\{0, a, b\}$, $J=\{0, c\}$. Then $I$ and $J$ are ideals of $M$. The zero-divisor graphs $\Gamma(M/I)$
	and $\Gamma(M/J)$, and
	the ideal-based zero-divisor graphs $\Gamma_{I}(M)$ and $\Gamma_{J}(M)$ are, respectively, drawed as follows:
	
	$$
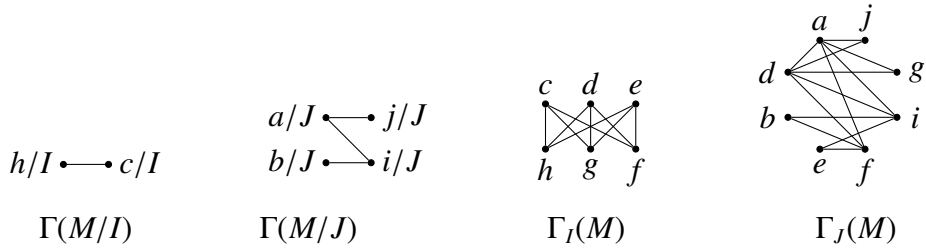
\begin{tikzpicture}
\tikzstyle{every node}=[draw,circle,fill=black,node distance=0.6cm,
minimum size=1.0pt, inner sep=0.8pt]
\node[circle] (1)                          [label=right : $c/I$]{};
\node[circle] (2) [left        of=1]       [label=left : $h/I$]{};	
\draw[-] (1) --   (2);
\end{tikzpicture}\quad\quad\quad
\begin{tikzpicture}
\tikzstyle{every node}=[draw,circle,fill=black,node distance=0.6cm,
minimum size=1.0pt, inner sep=0.8pt]
\node[circle] (1)                        [label=left :   $a/J$]{};
\node[circle] (2) [right        of=1]     [label=right : $j/J$]{};
\node[circle] (3) [below        of=1]     [label=left  : $b/J$] {};
\node[circle] (4) [below        of=2]     [label=right : $i/J$]{} ;

\draw[-] (1) --   (2); \draw[-] (1) --   (4);\draw[-] (3) --   (4);
\end{tikzpicture}\quad\quad\quad
\begin{tikzpicture}
\tikzstyle{every node}=[draw,circle,fill=black,node distance=0.6cm,
minimum size=1.0pt, inner sep=0.8pt]
\node[circle] (1)                          [label=above : $d$]{};
\node[circle] (2) [left        of=1]       [label=above : $c$]{};
\node[circle] (3) [right       of=1]      [label=above : $e$]{};
\node[circle] (4) [below        of=1]      [label=below  : $g$]{};
\node[circle] (5) [below        of=2]      [label=below : $h$]{};
\node[circle] (6) [below        of=3]      [label=below  : $f$]{};

\draw[-] (1) --   (4); \draw[-] (1) --   (5);\draw[-] (1) --   (6);
\draw[-] (2) --   (4); \draw[-] (2) --   (5); \draw[-] (2) --   (6);
\draw[-] (3) --   (4); \draw[-] (3) --   (5); \draw[-] (3) --   (6);
\end{tikzpicture}\quad\quad\quad
\begin{tikzpicture}
\tikzstyle{every node}=[draw,circle,fill=black,node distance=0.6cm,
minimum size=1.0pt, inner sep=0.8pt]
\node[circle] (1)                        [label=above :   $a$]{};
\node[circle] (2) [right        of=1]     [label=above : $j$]{};
\node[circle] (3) [below  left  of=1]     [label=left  : $d$] {};
\node[circle] (4) [below  right of=2]     [label=right : $g$]{} ;
\node[circle] (5) [below        of=3]     [label=left : $b$]{};
\node[circle] (6) [below        of=4]     [label=right : $i$] {};
\node[circle] (7) [below  right of=5]     [label=below : $e$]{} ;
\node[circle] (8) [below  left  of=6]     [label=below: $f$] {};

\draw[-] (1) --   (2); \draw[-] (1) --   (4);\draw[-] (1) --   (6);
\draw[-] (1) --   (8); \draw[-] (1) --   (3); \draw[-] (3) --   (2);
\draw[-] (3) --   (4); \draw[-] (3) --   (6); \draw[-] (3) --   (8);
\draw[-] (5) --   (6); \draw[-] (5) --   (8); \draw[-] (7) --   (6);
\draw[-] (7) --   (8);
\end{tikzpicture}
$$
$$\Gamma(M/I)\quad\quad\quad\quad\Gamma(M/J)\quad\quad\quad\quad\quad\quad\Gamma_{I}(M)\quad\quad\quad\quad\quad\quad \Gamma_{J}(M)$$

It is easy to see that $\Gamma_{I}(M)\cong K_{3, 3}$, $\diam(\Gamma(M/I))=1$,
$\diam(\Gamma_{I}(M))=2$ and $\diam(M/J)= \diam(\Gamma_{J}(M))=3$.

\end{example}

\subsection{Relationship between $\Gamma_{I}(A)$ and $\Gamma(A/I)$}\label{subsec:4.2}

In this subsection we will investigate some relationship between the ideal-based zero-divisor graph $\Gamma_{I}(A)$,
and the zero-divisor graph $\Gamma(A/I)$.
Let's first prove the following lemmas.
\begin{lemma}\label{lem:13}
		Let  $I\neq \{0\}$ be a proper ideal of an MV-algebra $A$ with $|A/I|\geq 4$.
	If $a, b\in V(\Gamma_{I}(A))$ such that $a/I\neq b/I$, then the following statements are equivalent:
	\begin{enumerate}
		\item   $a/I$ and $b/I$ are adjacent in $\Gamma(A/I)$.\label{it:4.11}
		\item
		$u$ and $v$ are adjacent in $\Gamma_{I}(A)$ for any
		$u\in a/I$ and any $v\in b/I$.
		\label{it:4.12}
		\item  $a$ and $b$ are adjacent in $\Gamma_{I}(A)$.\label{it:4.13}
	\end{enumerate}
\end{lemma}

\begin{proof}
	Suppose  $a, b\in V(\Gamma_{I}(A))$ such that $a/I\neq b/I$. It is clear that \mnoindent	
	\mref{it:4.12}$\Rightarrow $\mref{it:4.13}.
	
\mnoindent	
\mref{it:4.11}$\Rightarrow $\mref{it:4.12} Assume that \mref{it:4.11} holds. Then $a/I\odot b/I=0/I$. For any
 $u\in a/I$ and any $v\in b/I$, we have  $(u\odot v)/I=u/I\odot v/I=a/I\odot b/I=0/I$, and so $u\odot v\in I$.
	Consequently \mref{it:4.12} holds.

\mnoindent	
\mref{it:4.13}$\Rightarrow $\mref{it:4.11}	
	Assume that \mref{it:4.13} holds. Then $a\odot b\in I$, and so $a/I\odot b/I=(a\odot b)/I=0/I$. Hence \mref{it:4.11} holds.
\end{proof}

\begin{lemma}\label{lem:14} 	Let  $I\neq \{0\}$ be a proper ideal of an MV-algebra $A$ with $|A/I|\geq 3$,
  $a\in V(\Gamma_{I}(A))$ and $\Gamma_{I}(A)[a/I]$ be the subgraph of
	$\Gamma_{I}(A)$ induced by $a/I$. Then the following statements hold.
\begin{enumerate}
	\item  If $a/I\odot a/I=0/I$ in  $A/I$, then  $\Gamma_{I}(A)[a/I]$  is  complete.
	\label{it:4.21}
	
	\item	If $a/I\odot a/I\neq 0/I$ in  $A/I$, then  $\Gamma_{I}(A)[a/I]$  is  empty.
	\label{it:4.22}
\end{enumerate}
\end{lemma}
\begin{proof}
	Let $a\in V(\Gamma_{I}(A))$, and $x, y\in a/I$ with $x\neq y$. Then  $a/I\subseteq V(\Gamma_{I}(A))$ by Corollary \ref{cor:5}.
		
\mnoindent	
\mref{it:4.21}	If 
	$a/I\odot a/I=0/I$ in  $A/I$, then
	$(x\odot y)/I=x/I\odot y/I=a/I\odot a/I=0/I$,
	and so $x\odot y\in I$. Thus
	$\Gamma_{I}(A)[a/I]$ is complete.

\mnoindent	
\mref{it:4.22}	If $a/I\odot a/I\neq 0/I$ in  $A/I$, then
	$(x\odot y)/I=x/I\odot y/I=a/I\odot a/I\neq 0/I$,
	and so $x\odot y\not\in I$. Thus
	$\Gamma_{I}(A)[a/I]$ is empty.
\end{proof}

Recall that two simple graphs $G$ and $H$ are \name{isomorphic}, written $G\cong H$, if there is a
bijection $\theta :V(G) \rightarrow V(H)$ 
such that 
$u ~- ~v$ is an edge in $G$ if and only if $\theta(u) ~-~ \theta(v)$ is an edge in $H$.
Such a map $\theta$ is called \name{an isomorphism}
from $G$ to $H$ (see \cite[Exercise 1.2.5]{bon}).

\begin{theorem}\label{the:4.3}
	Let  $A$ and $B$ be finite MV-algebras such that $|A|=|B|$, and
	let  $I $ and $J$ be, respectively, proper ideals of  $A$ and $B$.
	If $A/I\cong B/J$, then $ \Gamma_{I}(A)\cong \Gamma_{J}(B)$.
\end{theorem}
\begin{proof}
	Let  $A$, $B$, $I$ and $J$ be as given. Assume that $A/I\cong B/J$. Then
	$\frac{|A|}{|I|}=|A/I|=|B/J|=\frac{|B|}{|J|}$ by Corollary \ref{cor:3.77},
	and consequently $|I|=|J|$.
	If $|I|=|J|=1$, then $I=\{0\}$ and $J=\{0\}$, which implies that $A\cong A/I\cong B/J\cong B$, and so
	$\Gamma_{I}(A)=\Gamma(A)\cong \Gamma (B)=\Gamma_{J}(B)$.	
	
	If  $|A/I|=|B/J|=2$, then $A/I\cong L_{2}\cong B/J$, and so  $ \Gamma_{I}(A)$ and $ \Gamma_{J}(B)$ are null  by Proposition \ref{pro:11}. Thus $\Gamma_{I}(A)\cong \Gamma_{J}(B)$.
	
	So we
	assume that $|A/I|=|B/J|=m+2$, where $m\geq 1$ is an integer.
	Put 
	\begin{center}
	 $A/I=\{0/I, a_{1}/I, a_{2}/I, \cdots, a_{m}/I, 1/I\}$ and
	$B/J=\{0/J, b_{1}/J,  b_{2}/J, \cdots, b_{m}/J, 1/J\}.$	
	\end{center}
Since  $A/I\cong B/J$, there exists an isomoprhism
	 $\theta: A/I\rightarrow B/J$  with
	$\theta(0/I)=0/J$, $\theta(1/I)=1/J$. Without loss of generality, we assume that
	$\theta(a_{i}/I)=b_{i}/J, i=1, 2, \cdots, m$.
	
	Next we construct a bijective map from $ V(\Gamma_{I}(A))$ to $ V(\Gamma_{J}(B))$.
	For any $a\in V(\Gamma_{I}(A))$ and any $b\in V(\Gamma_{J}(B))$, we have by Proposition \ref{pro:15} that
	$|a/I|=|I|=|J|= |b/J|$, and so
	there exists a bijective map
	$\eta_{i}: a_{i}/I\rightarrow b_{i}/J$ for any $ i\in \{1, 2, \cdots, m\}$.
	Notice that
	$V(\Gamma_{I}(A))=A\backslash (I\cup I^{*})=\bigcup_{i=1}^{m}a_{i}/I$ and  $V(\Gamma_{J}(B))=B\backslash (J\cup J^{*})=\bigcup_{i=1}^{m}b_{i}/J$.
	Define $$\eta=\bigcup_{i=i}^{m}\eta_{i}: V(\Gamma_{I}(A))\rightarrow V(\Gamma_{J}(B));\quad
	x\mapsto \eta_{i}(x),$$ when $x\in a_{i}/I, i=1, 2, \cdots, m$.
	Then it is easy to see that $\eta$ is a bijective map.
	Notice that $\theta: A/I\to B/J$ is an isomrophism.
	From Lemma \ref{lem:14} we have the induced subgraph $\Gamma_{I}(A)[a_{i}/I]$
	is complete (or empty) if and only if the induced subgraph
	$\Gamma_{J}(B)[b_{i}/J]$ is complete (or empty), which together by Lemma \ref{lem:13}, implies that $\eta$ is an isomorphism from
	$\Gamma_{I}(A)$ to $ \Gamma_{J}(B)$.
\end{proof}

If we replace finite by infinite,  Theorem \ref{the:4.3} doesn't necessarily hold (see Example \ref{ex:4.4}).

\begin{example}\label{ex:4.4}
	 Let $A=\mathcal{C}\times L_{2}$ and $B=\mathcal{C}\times L_{3}$, where  $\mathcal{C}$ is the MV-chain in Example \ref{exa:20}. Then it is clear that
	$|A|=|B|=|\mathcal{C}|$.
	Put $I=\{0\}\times L_{2}$
	and $J=\{0\}\times L_{3}$. Then $I$ and $J$ are, respectively,  ideals of $A$ and $B$. Also $A/I\cong \mathcal{C}\cong B/J$ by Corollary \ref{cor:6}.
	
	We will prove that $\Gamma_{I}(A)$ is not isomorphic to $\Gamma_{J}(B)$. Suppose, on the contrary, that
	$\Gamma_{I}(A)\cong\Gamma_{J}(B)$. Let
	$\delta$ be an isomorphism from 	$\Gamma_{I}(A)$ to $\Gamma_{J}(B)$.
	Then $\delta: V(\Gamma_{I}(A))\to V(\Gamma_{J}(B))$ is a bijective map.
	We claim that
	$\delta((c, 0)/I)=(c, 0)/J$, where
	$\delta((c, 0)/I)=\{\delta((x, y))~|~ (x, y)\in (c, 0)/I\}$.
	
	In fact, let $(a, z)\in  (c, 0)/I$. Then $(a, z)\in V(\Gamma_{I}(A))$ by Corollary \ref{cor:5}, and so $$\delta((a, z)) \in V(\Gamma_{J}(B))=B\backslash (J\cup J^{*}),$$  which implies that $\delta((a, z))\not \in J\cup J^{*}=\{0, 1\}\times L_{_{3}}$.
	So $\delta((a, z))=(b, y)$ for some $b\in C\backslash \{0, 1\}$ and $y\in L_{3}$.
	
	For any $(p, q)\in V(\Gamma_{I}(A))\backslash \{(a, z)\}$, 
		$(p, q)$ and $(a, z)$ are adjacent in $\Gamma_{I}(A)$. Indeed, by Proposition
		\ref{pro:11}, we have $(p, q)\not\in I\cup I^{*}=\{0, 1\}\times L_{2}$. Notice that $c\odot x=0$ for any $x\in \mathcal{C}\backslash \{0, 1\} $, we have $(c, 0)\odot (p, q)=(0, 0)$. 
		If $(a, z)/I=(p, q)/I$, then $(c, 0)/I=(p, q)/I$, and so	$(p, q)$ and $(a, z)$ are adjacent in $\Gamma_{I}(A)$
	by Lemma  \ref{lem:14} \ref{it:4.21}. If $(a, z)/I\neq (p, q)/I$, 
		then we have
		by Lemma  \ref{lem:13} that 
		$(p, q)$ and $(a, z)$ are adjacent in $\Gamma_{I}(A)$.
	 
 Since $\delta$ is an isomorphism from 	$\Gamma_{I}(A)$ to $\Gamma_{J}(B)$,
 it follows that
	$\delta((a, z))$ and $(u, v)$ are adjacent in $\Gamma_{J}(B)$
	for any $(u, v)\in V(\Gamma_{J}(B))\backslash \{\delta((a, z))\}$. 
	Thus, by Lemma \ref{lem:13},  $\delta((a, z))/J$ and $(u, v)/J$ are adjacent in  $\Gamma(B/J)$, which implies that  $\delta((a, z))/J \odot (u, v)/J=(0, 0)/J$, and so $\delta((a, z))\odot (u, v)\in J$. Recall that $\delta((a, z))=(b, y)$, where $b\in C\backslash \{0, 1\}$ and $y\in L_{3}$.
	If $b=c^{*}$, then $\delta((a, z))=(c^{*}, y)\neq ((2c)^{*}, 0)$, and so $((2c)^{*}, 0)\in  V(\Gamma_{J}(B))\backslash \{\delta((a, z))\}$, which implies that $\delta((a, z))\odot ((2c)^{*}, 0)\in J$. But
	 $$\delta((a, z))\odot ((2c)^{*}, 0)=(c^{*}, y)\odot ((2c)^{*}, 0)=((3c)^{*}, 0)\not\in\{0\}\times L_{3}=J,$$
	a contradiction. Thus $b\neq c^{*}$, and hence  $\delta((a, z))=(b, y)\neq (c^{*}, 0)$. It follows that $$(c^{*}, 0)\in  V(\Gamma_{J}(B))\backslash \{\delta((a, z))\},$$ and so $(b\odot c^{*}, 0)=(b, y)\odot (c^{*}, 0)=\delta((a, z))\odot (c^{*}, 0)\in J=\{0\}\times L_{3}$. 
Consequently $b\odot c^{*}=0$. 
	By Lemma \ref{lem:1} we get that
	$b\leq c$, and thus $b=c$ since $ b\in \mathcal{C}\backslash \{0, 1\}$ and
	$c$ is an atom in  MV-chain $\mathcal{C}$.
	Therefore $\delta((a, z))=(b, y)=(c, y)\in (c, 0)/J$, which implies that $\delta((c, 0)/I)\subseteq (c, 0)/J$.
	Similarly, we can prove that $\delta^{-1}((c, 0)/J)\subseteq (c, 0)/I$, i.e.
	$(c, 0)/J \subseteq\delta((c, 0)/I)$.
	Thus $\delta((c, 0)/I)= (c, 0)/J$, and hence  $|\delta((c, 0)/I)|= |(c, 0)/J|$.
	But  $|\delta((c, 0)/I)|=|(c, 0)/I|=|I|=2$ and $|(c, 0)/J|=|J|=3$  by Proposition $\ref{pro:15}$, a contradiction.
\end{example}

\subsection{Diameter of $\Gamma_{I}(A)$}\label{subsec:4.3}

In this subsection we will give some relationship between $\diam (\Gamma (A/I))$ and $\diam (\Gamma_{I}(A))$.

\begin{theorem}\label{the:2}	Let  $I$ be a proper ideal of an MV-algebra $A$ with $|A/I|\geq 3$.
	Then the following statements are equivalent:
	\begin{enumerate}
		\item $\diam(\Gamma_{I}(A))=0$.
		\label{it:4.51}
		\item  $\Gamma_{I} (A)$ is empty.
		\label{it:4.52}
		\item $|V(\Gamma_{I} (A))|=1$.
		\label{it:4.53}
		\item $I=\{0\}$  and $A\cong L_{3}$.
		\label{it:4.54}
		\item $I=\{0\}$  and  $\Gamma (A/I)$ is empty.
		\label{it:4.55}
	\end{enumerate}
\end{theorem}
\begin{proof}
	Since $\Gamma_{I}(A)$ is connected (see Theorem \ref{the:1}), we have
\ref{it:4.51}	$\Leftrightarrow $ \ref{it:4.52}$\Leftrightarrow$ \ref{it:4.53}.

\mnoindent			
\mref{it:4.54}$\Rightarrow $\mref{it:4.53} is clear.

\mnoindent	
\mref{it:4.53}$\Rightarrow $\mref{it:4.54}	
	Assume that  $|V(\Gamma_{I} (A))|=1$.
	Let $V(\Gamma_{I} (A))=\{a\}$. Then $a/I\subseteq V(\Gamma_{I} (A))$
	by Corollary \ref{cor:5}, and so
	$|a/I|=1$. It follows by Lemma \ref{lem:10} \ref{it:2.113}
	that $I=\{0\}$, and hence $\Gamma_{I} (A)=\Gamma(A)$. Therefore $A\cong L_{3}$
	by   \cite[Theorem 2]{gan},
	 \ref{it:4.54} holds.
	
\mnoindent	
\mref{it:4.54}$\Leftrightarrow $\mref{it:4.55}		follows  by  \cite[Theorem 2]{gan}.
\end{proof}

\begin{proposition}\label{pro:3333}
		Let $I$ be a proper ideal of an MV-algebra $A$. Then 
	$\diam (\Gamma (A/I))\leq \diam(\Gamma_{I}(A))$.
\end{proposition}
\begin{proof}
	If $\diam (\Gamma_{I} (A))=0$, then $\diam (\Gamma (A/I))=0$ by Theorem \ref{the:2}.
	
Assume that $\diam (\Gamma_{I} (A))=r\geq 1$.
	Let $a/I, b/I\in V(\Gamma (A/I))$ with $a/I \neq b/I$. Then $a\neq b$, and $a, b\in V(\Gamma_{I}(A))$ by Corollary \ref{cor:5}.
	So
	there exist $x_{1}, x_{2}, \cdots, x_{k}\in V(\Gamma_{I}(A))$
	such that $k\leq r-1$ and there is a path $a ~-~ x_{1} ~-~ x_{2} ~-~\cdots ~-~x_{k} ~-~b$ in $\Gamma_{I}(A)$.
	It follows  by Corollary \ref{cor:5} and Lemma \ref{lem:13}
	that there is a walk $a/I ~-~ x_{1}/I ~-~ x_{2}/I~-\cdots-~x_{k}/I ~-~b/I$ in $\Gamma(A/I)$, and
	consequently
	$\D_{\Gamma(A/I)}(a/I, b/I)\leq k+1\leq r$. Hence $\diam (\Gamma (A/I))\leq \diam(\Gamma_{I}(A))$.
\end{proof}

In Example \ref{exm:1000}, we have $\diam (\Gamma (M/I))=1$ and $ \diam(\Gamma_{I}(M))=2$.
So $\diam (\Gamma (A/I))$ may be strictly less than $ \diam(\Gamma_{I}(A))$ for an MV-algebra $A$.

\begin{theorem}\label{the:3}
	Let $\Gamma_{I}(A)$ be non-empty. Then the following statements are equivalent:
	\begin{enumerate}
		\item  $\diam(\Gamma_{I}(A))=1$;
		\label{it:4.71}
		\item   $\Gamma_{I} (A)$ is a complete graph;
		\label{it:4.72}	
		\item either $I\neq \{0\}$ and $A/I\cong L_{3}$;
		or  $I=\{0\}$ and $|A|=4$;
			\label{it:4.73}
		\item either $I\neq \{0\}$ and $\Gamma(A/I)$ is empty;
		or  $I=\{0\}$ and $\Gamma(A/I)\cong K_{2}$.
			\label{it:4.74}
	\end{enumerate}
\end{theorem}

\begin{proof}
	Suppose that  $\Gamma_{I}(A)$ is non-empty.
	It is clear that \ref{it:4.71}
	$\Rightarrow$ \ref{it:4.72}.
	
\mnoindent	
\mref{it:4.72}$\Rightarrow $\mref{it:4.73}	 Assume that \ref{it:4.72} holds. If
	$|A/I|\geq 5$, then by Lemma \ref{lem:5} and Proposition \ref{pro:11},
	there exist
	$a, b\in V(\Gamma_{I}(A))$ such that $a/I, a^{*}/I$ and $b/I$ are mutually distinct vertices in $\Gamma (A/I)$.
	Since $\Gamma_{I} (A)$ is a complete graph and $b^{*}/I$ is a vertex of $\Gamma(A/I)$, we have $a\odot b\in I$ and $a^{*}\odot b^{*}\in I$,
	and so $d(a^{*}, b)=(a\odot b)\oplus (a^{*}\odot b^{*})\in I$, which implies that $a^{*}/I=b/I$, a contradiction.
	Hence $|A/I|\leq 4$.
	
	If $|A/I|=3$, then  $A/I\cong L_{3}$ by  \cite[Theorem 2]{gan}. Since $\Gamma_{I}(A)$ is non-empty,
	we have $I\neq \{0\}$ by Theorem \ref{the:2}.
	
	If $|A/I|=4$, then we can assume that
	$A/I=\{0/I, a/I, a^{*}/I, 1/I\}$ by Lemma \ref{lem:5}.
	If
	$I\neq\{0\}$, then by Lemma \ref{lem:10} \ref{it:2.113},
	there exist $b\in a/I$ and $b^{*}\in a^{*}/I$ such that $a, b, a^{*}, b^{*}$ are mutually distinct
	vertices in $\Gamma_{I} (A)$.
	Since $\Gamma_{I} (A)$ is a complete graph, we have $a\odot b\in I$ and $a^{*}\odot b^{*}\in I$,
	so $d(a^{*}, b)=(a\odot b)\oplus (a^{*}\odot b^{*})\in I$, which implies $a^{*}/I=b/I=a/I$, a contradiction.
	Thus  $I=\{0\}$, and hence $|A|=|A/I|=4$.
	
\mnoindent	
\mref{it:4.73}$\Rightarrow $\mref{it:4.71}	
	If $I\neq \{0\}$ and $A/I\cong L_{3}$,
	say $A/I=\{0/I, a/I, 1/I\}$,
	then $V(\Gamma_{I}(A))=a/I$ by Proposition \ref{pro:11}, and
	$a/I=(a/I)^{\divideontimes}=a^{*}/I$. Thus  $\Gamma_{I} (A)$ is a complete graph
	by Lemma \ref{lem:14} \ref{it:4.21}. Since $I\neq \{0\}$, we have
	$|a/I|\geq 2$ by Lemma \ref{lem:10} \ref{it:2.113}, so  $\diam(\Gamma_{I}(A))=1$.
	
	If $I=\{0\}$ and $|A|=4$, then $\Gamma_{I}(A)=\Gamma(A)\cong K_{2}$, so
	$\diam(\Gamma_{I}(A))=1$.
	
\mnoindent	
\mref{it:4.73}$\Leftrightarrow$\mref{it:4.74} follows  by   \cite[Theorems 2 and 3]{gan}.
\end{proof}

 Proposition \ref{pro:22} describes the graph $\Gamma_{I} (A)$ when
$A$ is an MV-algebra with $|A/I|=4$.

\begin{proposition}\label{pro:22}
	Let $I$ be a proper ideal of an MV-algebra $A$ and $|I|=n\geq 2$.
	If $|A/I|=4$, then $\Gamma_{I}(A)\cong
	K_{n, n}$, or $\Gamma_{I}(A)\cong K_{n}\vee \emptyset_{n}$,
	the join  of  complete graph $K_{n}$ and  empty graph $\emptyset_{n}$.
\end{proposition}
\begin{proof}
	Let $I$ be a proper ideal of an  MV-algebra $A$ and $|I|=n\geq 2$.
	If $|A/I|=4$,  then $A/I\cong B_{4}$ or $L_{4}$ by  \cite[Theorem 3]{gan}.
	Assume that $A/I=\{0/I, a/I, a^{*}/I, 1/I\}$. We have $|a/I|=|a^{*}/I|=|I|=n$ by
	Proposition \ref{pro:15} and $V(\Gamma_{I}(A))=a/I \cup a^{*}/I$ by Proposition \ref{pro:11}.
	
	If $A/I\cong B_{4}$, then $a/I\odot a/I=a/I\neq 0/I$ and $a^{*}/I\odot a^{*}/I=a^{*}/I\neq 0/I$, so
	$\Gamma_{I}(A)\cong K_{n, n}$ by Lemmas \ref{lem:13} and  \ref{lem:14}.
	
	If $A/I\cong L_{4}$, then $A/I$ is an MV-chain. Without loss of generality,  assume
	that $a/I< a^{*}/I$. Then
	$a/I\odot a/I=0/I$
	and $a^{*}/I\odot a^{*}/I=a/I\neq 0/I$, so
	$\Gamma_{I}(A)\cong K_{n}\vee \emptyset_{n}$ by Lemmas \ref{lem:13} and  \ref{lem:14}.
\end{proof}

Now let's  turn our attention to the case when $\diam (\Gamma_{I}(A))=2$.
\begin{proposition}\label{pro:40}
	Let $I\neq \{0\}$  be a proper ideal of an  MV-algebra $A$  with $|A/I|=4$.
	Then $\diam (\Gamma_{I}(A))=2$.
\end{proposition}
\begin{proof}
	Let $I, A$  be as given, say
	$A/I=\{0/I, a/I, a^{*}/I, 1/I\}$.
	Then $V(\Gamma_{I}(A))=a/I\cup a^{*}/I$ by Proposition \ref{pro:11}.
	To prove that $\diam (\Gamma_{I}(A))= 2$, let  $x, y\in V(\Gamma_{I}(A))$ with $x\neq y$. Consider the
	following cases:
	
Case (i).  If $x\in a/I$ and $y\in a^{*}/I$, or $y\in a/I$ and $x\in a^{*}/I$, then
	$x$ and $y$ are adjacent in $\Gamma_{I}(A)$ by Lemma \ref{lem:13}.
	
Case (ii). If $x, y\in a/I$, then there is a path $x~-~ a^{*}~-~ y$ in $\Gamma_{I}(A)$ by Lemma \ref{lem:13}.

Case (iii). if $x, y\in a^{*}/I$, then there is a path $x~-~ a ~- ~y$ in $\Gamma_{I}(A)$ by Lemma \ref{lem:13}.

	Thus we have $\D_{\Gamma_{I}(A)}(x, y)\leq 2$, and hence $\diam (\Gamma_{I}(A))\leq 2$. But $\diam (\Gamma_{I}(A))\neq 0$
	by Theorem \ref{the:2} and $\diam (\Gamma_{I}(A))\neq 1$ by Theorem  \ref{the:3},  so $\diam (\Gamma_{I}(A))=2$.
\end{proof}

\begin{proposition}\label{pro:23}
Let  $I\neq \{0\}$ be a proper ideal of an  MV-algebra $A$. If 
	 $\diam (\Gamma (A/I))=2$,
	then $\diam (\Gamma_{I}(A))=2$.
\end{proposition}
\begin{proof}
	Suppose that $I\neq \{0\}$ and
	$\diam (\Gamma (A/I))=2$. Then 
	 $|A/I|\geq 5$, sine $|A/I|\in \{3, 4\}$ implies that 	$\diam (\Gamma (A/I))=0$ or $1$. Also, by Proposition \ref{pro:3333}, we have  $2=\diam (\Gamma (A/I)) \leq \diam (\Gamma_{I}(A))$.
	Let $x, y\in V(\Gamma_{I}(A))$ such that $x\neq y$ and there is no edge between $x$ and $y$ in $\Gamma_{I}(A)$. 
	
	If $x/I=y/I$, then there exists $ z/I\in V(\Gamma(A/I))$ such that
	$z/I\neq x/I$, and
	$x/I$ and $ z/I $ are adjacent in $\Gamma(A/I)$, since $|V(\Gamma (A/I))|\geq 3$ and $\Gamma(A/I)$ is connected.
	It follows by Lemma \ref{lem:13} that there is a path $x~-~z ~- ~y$ in  $\Gamma_{I}(A)$.
	
	If $x/I\neq y/I$, 
	then  there is no edge between $x/I$ and $y/I$ in $\Gamma(A/I)$ by Lemma \ref{lem:13}. Thus
	there exists $ a/I\in V(\Gamma(A/I))$ such that there is a path
	$x/I~ -~a/I~-~ y/I$ in $\Gamma(A/I)$, since $\diam (\Gamma (A/I))=2$. It follows by Lemma \ref{lem:13}
	that there is a path
	$x~-~a~-~y$ in  $\Gamma_{I}(A)$.
	
	From the above arguments, we get 
 $\diam (\Gamma_{I}(A))=2$.
\end{proof}

\begin{theorem}\label{the:5}
Let  $I\neq \{0\}$ be a proper ideal of an  MV-algebra $A$.
	Then $\diam (\Gamma_{I}(A))=2$ if and only if $\diam (\Gamma (A/I))=1$ or $2$.
\end{theorem}
\begin{proof}
	Let  $I\neq \{0\}$.
	If $\diam (\Gamma_{I}(A))=2$, then  $\diam (\Gamma (A/I))\leq 2$ by Proposition \ref{pro:3333}.
	Also, by Theorem \ref{the:2} we have $\diam (\Gamma (A/I))\neq 0$, so
	 $\diam (\Gamma (A/I))=1$ or $2$.
	
	Conversely, if $\diam (\Gamma (A/I))=1$, then  $|A/I|=4$ by  \cite[Theorem 3]{gan}, so
	$\diam (\Gamma_{I}(A))=2$ by Proposition \ref{pro:40}.
	If $\diam (\Gamma (A/I))=2$, then  $\diam (\Gamma_{I}(A))=2$ by Proposition \ref{pro:23}.
\end{proof}

\begin{theorem}\label{the:6} Let  $I$ be a proper ideal of an  MV-algebra $A$. Then
	$\diam (\Gamma_{I}(A))=3$ if and only if $ \diam (\Gamma (A/I))=3$.
\end{theorem}
\begin{proof}
	If $I=\{0\}$, then $A/I\cong A$ and
	$\Gamma_{I}(A)\cong \Gamma(A)\cong \Gamma(A/I)$, so the result is true.
	
	Assume that $I\neq \{0\}$. If $\diam (\Gamma (A/I))=3$, then $\diam (\Gamma_{I}(A))=3$ by Theorem \ref{the:1} and Proposition
	\ref{pro:3333}. Conversely,
	if  $\diam (\Gamma_{I}(A))=3$, then
	$\diam (\Gamma (A/I))=3$ by Theorem \ref{the:3} and Theorem \ref{the:5}.
\end{proof}

\subsection{Girths of  $\Gamma_{I}(A)$}\label{subsec:4.4}
Let us first discuss the relationship between the girth of the zero-divisor graph $\Gamma(A/I)$
and the ideal-based zero-divisor graph of $\Gamma_{I}(A)$.
\begin{proposition}\label{pro:5.18}
	Let $I$ be a proper ideal of an MV-algebra $A$. Then 	
	$\gr (\Gamma_{I} (A))\leq \gr(\Gamma(A/I))$.
	
	 In particular, if 	
	$\gr(\Gamma(A/I))=3$, then 	$\gr (\Gamma_{I} (A))=3$.
\end{proposition}
\begin{proof}
	If $\gr (\Gamma(A/I))=\infty$, then it is clear 	$\gr (\Gamma_{I} (A))\leq \gr(\Gamma(A/I))$.
		Assume that $\gr (\Gamma (A/I))=k$, where $k\geq 3$ is an integer, and let $a/I ~-~ x_{1}/I ~-~ x_{2}/I~- \cdots - ~x_{k-1}/I ~- ~a/I$ be a cycle of length $k$ in $\Gamma(A/I)$. Then
	$a ~-~ x_{1} ~-~ x_{2} ~-~\cdots ~-~x_{k-1} ~-~a$ is a cycle of length $k$ in $\Gamma_{I}(A)$ 	 by Corollary \ref{cor:5} and Lemma \ref{lem:13}. Hence 	$\gr (\Gamma_{I} (A))\leq gr(\Gamma(A/I))$.
	
	In particular, if $\gr (\Gamma (A/I))=3$,  and let $a/I ~-~ x_{1}/I ~-~ x_{2}/I ~- ~a/I$ be a triangle in $\Gamma(A/I)$. Then
	$a ~-~ x_{1} ~-~ x_{2} ~-a$ is a triangle in $\Gamma_{I}(A)$ 	 by Corollary \ref{cor:5} and Lemma \ref{lem:13}. Hence 	$\gr (\Gamma_{I} (A))=3$.
\end{proof}

In Example \ref{exm:1000}, we have $\gr (\Gamma (M/I))=\gr(\Gamma (M/J))=\infty$, while  $ \gr(\Gamma_{I}(M))=4$ and  $ \gr(\Gamma_{J}(M))=3$.
So $\gr (\Gamma (A/I))$ may be strictly large than $ \gr(\Gamma_{I}(A))$ for an MV-algebra $A$.

\begin{lemma}\label{lem:15}
	Let $I$ be a proper ideal of an MV-algebra $A$ and $I\neq \{0\}$.	If there exists $a/I\in A/I \backslash \{0/I, 1/I\}$ such that $a/I< a^{*}/I$ in $A/I$, then 	
	$\gr (\Gamma_{I} (A))=3$. 
	
	In pariticular, if
	$A/I$ is an MV-chain with $|A/I|\geq 4$, then
		$\gr (\Gamma_{I} (A))=3$.
\end{lemma}
\begin{proof}
Let $I$ and $a/I$ be as given. Then
	$a/I\odot a/I=0/I$ in  $A/I$ by Lemma \ref{lem:1}, and so
	$\Gamma_{I}(A)[a/I]$ is complete by Lemma \ref{lem:14} \ref{it:4.71}. Thus 	$\gr (\Gamma_{I} (A))=3$
	when $|I|\geq 3$. 
	
	If $|I|=2$, then $|a/I|=|a^{*}/I|=|I|=2$ by Proposition \ref{pro:15}. Fix
	$b\in a/I$ and $b\neq a$. We get a triangle
	$b~-~a ~-~a^{*}~-~b$ by Lemma \ref{lem:13} and Lemma \ref{lem:14}. Hence $\gr (\Gamma_{I} (A))=3$.
\end{proof}

\begin{theorem}\label{the:5.20}
	Let $I$ be a proper ideal of an MV-algebra $A$, $I\neq \{0\}$ and  $|A/I|\geq 3$. Then
	$$\gr (\Gamma_{I}(A))=
	\begin{cases}
	\infty,  & \textrm{if}~ A\cong L_{2}\times L_{3} ~ and~~ I\cong L_{2}\times \{0\}; \\
	4,  & \textrm{if}~ A/I\cong B_{4}; \\
	
	3,  & \textrm{otherwise}.
	\end{cases}$$
\end{theorem}
\begin{proof}
	If  $A/I\cong L_{3}$, say $A/I=\{0/I, a/I, 1/I\}$,  then 
	  $\Gamma_{I}(A)$ is complete by Corollary \ref{cor:5} and Lemma \ref{lem:14} \ref{it:4.21}, and so
	$\gr(\Gamma_{I}(A))=\infty$ when $|I|=2$; and
	$\gr(\Gamma_{I}(A))=3$ when $|I|\geq 3$. Moreover, when $|A/I|=3$ and $|I|=2$, we have
	$|A|=|A/I||I|=6$ by Corollary \ref{cor:3.77}, so $A$ is isomorphic to $ L_{6}$ or $L_{2}\times L_{3}$ by \cite[Theorem 7]{gan}.
	If $A\cong L_{6}$, then $A$ does not contain an ideal such that $|I|=2$, since $L_{6}$ is simple by Lemma \ref{lem:44}. If $A\cong L_{2}\times L_{3}$ and $|I|=2$, then $I\cong L_{2}\times \{0\}$  by  Lemma \ref{lem:44} and Proposition \ref{cor:10000}.
	
 If $A/I\cong L_{4}$, then  $\gr(\Gamma_{I}(A))=3$ by Lemma \ref{lem:15}.
	 If $A/I\cong B_{4}$, then $\Gamma_{I}(A)$ is a complete biparitite graph by Lemma \ref{lem:13}
	 and Lemma \ref{lem:14}, so $\gr(\Gamma_{I}(A))=4$.
	 	
	If  $A/I\cong L_{5}$, then $\gr(\Gamma_{I}(A))=3$ by Lemma \ref{lem:15}.

	If  $A/I\cong L_{2}\times L_{3}$,   then there exist $a/I\in A/I \backslash \{0/I, 1/I\}$ such that $a/I< a^{*}/I$ in $A/I$,
	since $(0, \frac{1}{2})<(1, \frac{1}{2})=(0, \frac{1}{2})^{*}$ in $L_{2}\times L_{3}$. Thus	
	$\gr (\Gamma_{I} (A))=3$ by Lemma \ref{lem:15}. 
	
For other cases, by Theorem \ref{the:4.30} and Proposition \ref{pro:5.18}, we have 	$\gr (\Gamma_{I} (A))=3$. 
\end{proof}

\vspace{3mm}

\noindent
{\bf Acknowledgments.}
This work is supported by the NSFC Grants (Nos. 12171022, 12261001) and the Guangxi Natural Science Foundation (Grant No. 2021GXNSFAA220043) and High-level talents for scientific research of Beibu Gulf University (2020KYQD07). We thanks the anonymous referees for helpful their suggestions.

\end{document}